\documentclass[12pt]{elsarticle}

\usepackage{graphicx}
\usepackage{amsmath}
\usepackage{amsthm}
\usepackage{amsfonts}
\usepackage[brazil]{babel}
\usepackage[utf8]{inputenc}
\usepackage{graphicx}
\usepackage{amsmath}
\usepackage{amsthm}
\usepackage{amsfonts}
\usepackage[utf8]{inputenc}

\newfont{\Mb}{msbm10}
\newcommand{\C}{\mbox{\Mb\symbol{67}}}
\newcommand{\R}{\mbox{\Mb\symbol{82}}}

\newcounter{bla}

\newtheorem{teor}{Theorem}[section]
\newtheorem{cor}{Corollary}[section]
\newtheorem{obs}{Remark}[section]
\newtheorem{defin}{Definition}[section]

\newtheorem{algor}{Algorithm}[section]

\journal{Journal of Computational Physics}

\begin{document}

\title{\uppercase{Generalized Liénard equations presenting Liouvillian first integrals}}

\author[uerj]{L.G.S. Duarte}
\ead{lgsduarte@gmail.com.br}

\author[uerj]{L.A.C.P. da Mota\corref{cor1}}
\ead{lacpdamota@uerj.br or lacpdamota@gmail.com}

\author[uerj,imetro]{A.B.M.M.  Queiroz}
\ead{andream.melo@gmail.com}

\cortext[cor1]{Corresponding author {\footnotesize \newline L.G.S. Duarte and L.A.C.P. da Mota wish to thank Funda\c c\~ao de Amparo \`a Pesquisa do Estado do Rio de Janeiro (FAPERJ) for a Research Grant.}}

\address[uerj]{Universidade do Estado do Rio de Janeiro,
{\it Instituto de F\'{\i}sica, Depto. de F\'{\i}sica Te\'orica},
{\it 20559-900 Rio de Janeiro -- RJ, Brazil}}

\address[imetro]{Inmetro - Instituto Nacional de Metrologia, Qualidade e Tecnologia.}

\begin{abstract}
In this work we propose a new method to determine Liouvillian first integrals for Liénard–Levinson–Smith (LLS) equations. Our procedure presents three interesting characteristics: In the case where the equation depends on parameters, the method allows an analysis of the integrability region of the equation's parameters in a natural way; moreover, as the Liouvillian first integrals are directly linked to the existence of Darboux integrating factors, the procedure allows an analysis of the equation's algebraic limit cycles (i.e., its periodic algebraic solutions in the fase space); finally, as the LLS (not forced) equations can be transformed into general first-order ordinary differential equations (1ODEs), the method provides (as a bonus) an efficient way to solve 1ODEs with elementary functions that present a Liouvillian general solution. The method is implemented in a Maple package {\it LeapsLS} and it includes commands that allow obtaining all the intermediate steps of the process of finding the Liouvillian first integrals.
\end{abstract}


\maketitle

\newpage

\section{Introduction}
\label{introd}

\medskip

The Liénard-type differential equation \cite{Liena} 
\begin{equation}
\label{lien}
\ddot{x} + F(x)\,\dot{x} + G(x) = 0.
\end{equation}
and its generalization of the form (the so called Liénard-Levinson-Smith (LLS) equations) \cite{LeSm}
\begin{equation}
\label{levsmi}
\ddot{x} + F(x,\dot{x})\,\dot{x} + G(x) = 0,
\end{equation}
appear in several areas of knowledge such as physics, chemistry, biology, engineering etc (see \cite{SaGaRa,ChSeLa,Val} and references therein)  and, therefore, are considered as one of the most important classes of second-order nonlinear ordinary differential equations (nonlinear 2ODEs)\footnote{These equations model a wide range of phenomena presenting nonlinear behavior and, probably, the first famous Liénard-type differential equation was obtained by van der Pol \cite{Pol,Pol2,Pol3}, who was studying nonlinear oscillations in electrical circuits. For a derivation of the van der Pol equation see \cite{Gin}.}. The LLS equations also appears when one is performing traveling wave reduction of nonlinear partial differential equations (examples include the Burgers–Korteweg–de Vries and Fisher equations -- see \cite{KuZa,GaBrRo} and references therein). 
If we make $y=\dot{x}$, we can write the LLS equation as an autonomous system (usually called Liénard differential system or LLS system)
\begin{equation}
\label{llssys}
\left\{ \begin{array}{l}
\dot{x} = y, \\
\dot{y} = -\left(F(x,y)\,y + G(x)\right).
\end{array} \right.
\end{equation}
The most important problems of research in the Liénard differential systems consist of the study of their limit cycles (i.e., its periodic orbits) and the integrability problem (see \cite{ChrL,LlMeTe,PaBiSeLa,LlVaL,ZhWa,GiLl,Gin,LlZh,Dem,CaGu} and references therein). Therefore, if $F(x,y)$ and $G(x)$ are rational functions of $(x,y)$, analyzing the existence of Liouvillian first integrals of LLS system (\ref{llssys}), we will be reaching a double objective, since the presence of a Liouvillian first integral implies the existence of a Darboux integrating factor of the type \cite{Chr}
\begin{equation}
\label{LiouvIF}
R = \mathrm{\large e}^{Z_0} \prod_i {p_i}^{n_i}
\end{equation}
where $Z_0$ is a rational function of $(x,y)$, the $p_i$ are irreducible polynomials of $(x,y)$ and the $n_i$ are constants. Furthermore, the denominator of $Z_0$ and the $p_i$ are Darboux polynomials of the LLS system (\ref{llssys}) and we can use them (the Darboux polynomials) to look for algebraic limit cycles of the system (periodic solutions). Besides that, if the functions $F(x,y)$ and $G(x)$ present an elementary function $\theta(x,y)$ we can build a strategy (construct an algorithm) to find a Liouvillian first integral for (\ref{llssys}).

Regarding the search for Liouvillian first integrals\footnote{The methods that use the Lie symmetry approach also have their place. See, for example \cite{Olv,Ibr,BlAn,Sch,Noscpc1997,Noscpc1998,AbGu,GoLe,AdMa,GaBrSe,MuRo,MuRo2,PuSa,Nuc}.}, the methods that use the Darboux-Prelle-Singer approach stand out. In this `direction' we can highlight: \cite{Dar,PrSi,Sin,Chr,CaLl,ChLlPaZh,Nosjpa2002-1,Nosjpa2002-2,Noscpc2002,ChGiGiLl,Nosjcam2005,Noscpc2007,ChLlPe,ChLlPaWa2,ChLlPaWa3,Zha,FeGi,
BoChClWe,ChCo,Nosjde2021}.
We have been studying the existence of (and methods to search for) Liouvillian first integrals of vector fields in $\R^2$ since 2002 \cite{Nosjpa2002-1,Nosjpa2002-2,Noscpc2002,Nosjcam2005,Noscpc2007,Nosjde2021}. In this paper, we combine the technique we used in \cite{Noscpc2019}\footnote{In that work we use the (so-called) $S$-function (introduced in \cite{Nosjpa2001} and used by us and other people \cite{Noscpc2019,Nosamc2007,LaRa,ChSeLa,Nosjmp2009,Nosjpa2010,GoPiSe,TPCSL,MCSL,tito1,tito2}) that improves greatly the efficiency on the searching for Liouvillian first integrals of 2ODEs.} with a reasoning very close to the idea used in \cite{Noscpc2002}\footnote{The central idea is to `transform' elementary functions into rational functions.}. Our method basically starts by assuming the existence of a Liouvillian first integral (of a certain type) for a LLS equation and comes to a conclusion: the existence (or not) of a Liouvillian first integral (of a certain type) for a polynomial vector field in three variables and, in the case of a positive answer, the determination of the LLS's first integral. Briefly, from the fact that the LLS equation (\ref{levsmi}) has the Lie symmetry $\partial_t$, we can make the transformation
\begin{equation}
\label{tlevsmi}
\left\{ \begin{array}{l}
x = x, \\
\dot{x} = y, \\
\ddot{x}= \frac{dy}{dx} \, y,
\end{array} \right.
\end{equation}
and reduce its order, i.e., we apply the transformation (\ref{tlevsmi}) to the LLS equation (\ref{levsmi}) leading to the 1ODE given by
\begin{equation}
\label{1levsmi}
\frac{dy}{dx}= - \frac{F(x,y)\,y + G(x)}{y}.
\end{equation}  
Now, if $ F $ presents an elementary function $\theta(x,y)$ (or if $ F $ and $ G $ present an elementary function $\theta(x)$), we make another transformation $\,\theta \,\rightarrow\,z\,$ that `leads' the 1ODE (\ref{1levsmi}) into a polynomial vector field $\cal X$ in three variables. If we find a Liouvillian first integral for $\cal X$ then we only have to perform the inverse transformation $\,z\,\rightarrow\,\theta \,$ to obtain the general solution of the 1ODE (\ref{1levsmi}) in implicit form. Finally, applying the inverse of the transformation (\ref{tlevsmi}), we obtain the desired Liouvillian first integral.

\bigskip

\noindent
This paper is organized as follows:

\medskip

\noindent
In the first section, we present some basic concepts involved in the Darboux-Prelle-Singer (DPS) approach and in the $S$-function method. 

\medskip

\noindent
In the second section, we show how we can associate a LLS equation (\ref{levsmi}) with an elementary function to a polynomial vector field in three variables and develop a technique similar to the one we built in \cite{Noscpc2019} (to search for Liouvillian first integrals of rational 2ODEs using the $S$-function method) in order to search for a Liouvillian first integral of the vector field. We use these two steps to propose a procedure that can succeed even in cases where the integrating factors have Darboux polynomials of very high degree. In the end of the section (after a few examples in order to clarify the steps of the procedure), we propose a semi algorithm to determine a Liouvillian first integral of a LLS equation presenting an elementary function and we describe it in a more formal way.

\medskip

\noindent
In the third section, we discuss the performance of the method\footnote{We made an implementation of the method in Maple -- see appendix}. We begin by presenting a set of LLS equations that our method can deal without problems but is not easily solved by the `traditional' methods. Finally, we show that, for an LLS equation presenting parameters, the method employed can perform an analysis of the integrability region of the parameters: We present an example that the package can successfully handle.

\medskip

\noindent
Finally, we present our conclusions.

\medskip

\noindent
In the appendix, we present a computacional package (in Maple) to search for Liouvillian solutions of LLS equations such that the functions $F(x,y)$ and $G(x)$ may present an elementary function $\theta(x,y)$ (i.e., $F(x,y,\theta(x,y))$ ; $G(x)$ ) or $\theta(x)$ (i.e., $F(x,y,\theta(x))$ ; $G(x,\theta(x))$). We start by giving an overview of the package and presenting a summary of the commands. Then, we go deeper into the description of each command and, in the last subsection, we present an example of the usage of the commands to show how they work in practice.

\section{Basic concepts and results}
\label{bcr}

In this section we will set the mathematical basis for presenting our method. In the first subsection we will present the Darboux-Prelle-Singer approach in a nutshell and, in the second, the basics of the $S$-function method.

\subsection{The Darboux-Prelle-Singer approach}
\label{dpsa}

When we talk succinctly about the Darboux-Prelle-Singer approach to seeking Liouvillian first integrals of autonomous systems of polynomial 1ODEs in the plane, we can summarize the central idea in one paragraph:

{\it Find the Darboux polynomials (DPs) associated with the 1ODE system and use them to determine an integrating factor.}

Without exaggeration, the determination of the DPs is the most important and (unfortunately) the most complex part (by far) of the whole process. Thus, with a couple of definitions and results, we can describe (briefly) the DPS approach when applied to vector fields in the plane. To begin, let's define more formally the concepts we referred to.
Consider the following system of polynomial 1ODEs in the plane:

\begin{equation}
\label{sys1}
\left\{ \begin{array}{l}
\dot{x} = f(x,y)\\ [2mm]
\dot{y} = g(x,y),
\end{array} \right.
\end{equation}
where $f$ and $g$ are coprime polynomials in $(x,y)$ and the dot means derivative with respect to a parameter $t$ ($\dot{u} \equiv du/dt$).

\begin{defin}
A function $I(x,y)$ is called a {\bf first integral} of the system {\em (\ref{sys1})} if $I$ is constant over the solutions of {\em (\ref{sys1})}.
\end{defin}

\begin{defin}
The vector field associated with the system {\em (\ref{sys1})}, 
\begin{equation}
\label{vf1}
X \equiv f(x,y)\,\partial_x+g(x,y)\,\partial_y,
\end{equation}
is called {\bf Darboux operator} associated with {\em (\ref{sys1})}.
\end{defin}

\begin{obs}
If $I(x,y)$ is a {\bf first integral} of the system {\em (\ref{sys1})}, then $X(I) = 0$.
\end{obs}

\begin{defin}
The polynomial $p(x,y) \in \C[x,y]$ is called a {\bf Darboux polynomial} of the vector field $X$ if $X(p) = q\,p$, where $q$ is a polynomial denominated {\bf cofactor}.
\end{defin}

\begin{teor}[Prelle-Singer]
\label{ps}
If the system {\em (\ref{sys1})} presents an elementary first integral $I$, then there exists an integrating factor $R$ for the system {\em (\ref{sys1})} of the form $R = \prod_i {p_i}^{n_i}$, where the $p_i$ are irreducible Darboux polynomials of system {\em (\ref{sys1})} and $n_i$ are rational numbers.
\end{teor}

\noindent
{\it Proof.} For a proof see \cite{PrSi}.

\bigskip

Since $\frac{X(R)}{R} = - \,{\rm {\bf div}}(X) = - (f_x + g_y)$, substituting $R = \prod_i {p_i}^{n_i}$ we obtain (see \cite{Noscpc2002})
\begin{equation}
\label{dpispi}
\sum_i \,n_i\,\frac{X(p_i)}{p_i} = \sum_i \,n_i\,q_i = - (f_x + g_y),
\end{equation}
where the polynomials $q_i$ ($\equiv \frac{X(p_i)}{p_i}$) are called {\em cofactors} of the Darboux polynomials $p_i$. Therefore, a possible method to find an elementary first integral is:

\bigskip

\noindent
{\bf Prelle-Singer method:}(sketch)
\begin{itemize}
\item Determine the DPs $p_i$ associated with the system.
\item Find numbers $n_i$ that satisfy $\sum_i \,n_i\,q_i = - (f_x + g_y)$.
\item Construct the integrating factor $R = \prod_i {p_i}^{n_i}$ and find a first integral $I(x,y)$ by quadratures.
\end{itemize}

\begin{obs}
If the system {\em (\ref{sys1})} presents a non elementary Liouvillian first integral the Prelle-Singer (PS) method needs some modifications (see the Christopher-Singer (CS) method in {\em \cite{Nosjcam2005,Noscpc2007,Nosjde2021}}). However, the main premisse remains the same: determine the DPs of the system.
\end{obs}

\begin{obs}
For vector fields in $\R^3$, things get a little more complicated: first, the integrating factor $R$ is no longer a Jacobi multiplier (as it happens in $\R^2$). Thus, $X (R) / R$ is no longer a known polynomial and depends on a rational function not determined a priori. Second, we need to determine Darboux polynomials in three variables (a very difficult task since it was no longer easy to determine them in two variables).
\end{obs}

\medskip

\subsection{The $S$-function method}
\label{sfm}

The concept of S-function emerged when we were trying to adapt the DPS approach to rational second order ordinary differential equations (rational 2ODEs). The idea was to `add' a null 1-form to the 1-form that represented the 2ODE (see \cite{Nosjpa2001,Nosjmp2009}). However, the problem of calculating DPs was so costly (in computational terms) that we developed a mixed procedure we called `the $S$-function method'. Briefly, the method did not use the computation of Darboux polynomials (in exchange for the resolution of two 1ODEs). The good news is that the $S$-function method proved to be more effective precisely in problems in which the Darboux polynomials had a relatively high degree, being a great alternative to be used in these cases. In this section, we will see what the method consists of and its basic operation.

\subsubsection{Some basic definitions and results}
\label{sbdr}

Consider the rational 2ODE given by
\begin{equation}
\label{2oder1}
z'=\frac{dz}{dx}=\phi(x,y,z)=\frac{M(x,y,z)}{N(x,y,z)},  \,\,(z \equiv y'),
\end{equation}
where $M$ and $N$ are coprime polynomials in $\C[x,y,z]$.

\begin{defin}
A function $I(x,y,z)$ is called {\bf first integral} of the 2ODE {\em (\ref{2oder1})} if $I$ is constant over the solutions of {\em (\ref{2oder1})}.
\end{defin}

\begin{obs} If $I(x,y,z)$ is a first integral of the 2ODE {\em (\ref{2oder1})} then, over the solution curves of {\em (\ref{2oder1})}, the exact 1-form $\omega:=dI=I_x\,dx+I_y\,dy+I_{z}\,dz\,$ is null.
\end{obs}
\noindent
Over the solution curves of (\ref{2oder1}), we have two independent null 1-forms:
\begin{eqnarray}
\label{alfa2}
\alpha & := & \phi\,dx-dz, \\
\label{beta2}
\beta & := & z\,dx-dy.
\end{eqnarray}
So, the 1-form $\omega$ is in the vector space sppaned by the 1-forms $\alpha$ and $\beta$, i.e., we can write $\omega = r_1 (x,y,z)\,\alpha+r_2 (x,y,z)\,\beta$:
\begin{eqnarray}
\label{eqform2edo1}
dI &=& I_x\,dx+I_y\,dy+I_{z}\,dy = r_1\,(\phi\,dx - dz)\,+r_2\,(z\,dx - dy) \nonumber \\
    &=& (r_1\,\phi\,+r_2\,z)dx + (-r_2)dy + (-r_1)dz.
\end{eqnarray}
\noindent
If we call $\,r_1 \equiv R\,$ and $\,\frac{r_2}{r_1} \equiv S\,$, we can write
\begin{equation}
\label{dIrands}
dI = R\,\left[(\phi\,+z\,S)dx - S\,dy - dz\right].
\end{equation}

\noindent
Therefore, we have that $I_x = R\,(\phi\,+z\,S), \, I_y = -R\,S, \,I_z = -R$. 

\begin{defin}
\label{fatint1form}
Let $\gamma$ be a 1-form. We say that $R$ is an {\bf integrating factor} for the 1-form $\gamma$ if $R\,\gamma$ is an exact 1-form.
\end{defin}

\begin{defin}
Let $I$ be a first integral of the {\em 2ODE (\ref{2oder1})}. The function defined by $S := I_y/I_z$ is called a \mbox{\boldmath $S$}{\bf -function} associated with the {\em 2ODE} through the first integral $I$.
\end{defin}

\begin{obs}
From the above definitions and in view of {\em (\ref{dIrands})} we can see that $R$ is an integrating factor for the 1-form $\,(\phi\,+z\,S)dx - S\,dy - dz\,$ and $S$ is a $S$-function associated with the {\em 2ODE (\ref{2oder1})} through $I$ .
\end{obs}

\begin{teor}
\label{concom}
Let $\,I(x,y,z)\,$ be a first integral of the {\em 2ODE (\ref{2oder1})}. If $S$ and $R$ are as in {\em (\ref{dIrands})}, then we can write
\begin{eqnarray}
\label{eqproof1r}
&& D_x(R) +R\,(S+\phi_z)=0\,, \\
\label{eqproof2r}
&& S\,D_x(R) + R\,( D_x(S)+\phi_y)=0\,, \\
\label{eqproof3r}
&& -(R_z\,S+R\,S_z) + R_y=0\,,
\end{eqnarray}
where $\,D_x := \partial_x + z\,\partial_y + \phi(x,y,z)\,\partial_z\,$.
\end{teor}

\noindent
{\it Proof:} See \cite{Nosjmp2009}.

\begin{cor}
\label{corS}
Let $S$ be a $S$-function associated with the {\em 2ODE} $\,z'=\phi(x,y,z)$. Then $S$ obeys the following equation:
\begin{equation}
\label{eqS}
D_x(S)=S^2+\phi_{z}\,S-\phi_y\,.
\end{equation}
\end{cor}

\noindent
{\it Proof:} See \cite{Noscpc2019}

\subsubsection{The Associated 1ODEs}
\label{conectSandS}
From (\ref{eqS}) we can see that a $S$-function associated with the rational 2ODE (\ref{2oder1}) satisfies a quasilinear 1PDE in the variables $(x,y,z)$. 
Over the solutions of the 2ODE (\ref{2oder1}) we have that $y=y(x)$ and $z=z(x)$ and, therefore, the operator $D_x$ is, formally, ${\frac{d}{dx}}$. So, formally, over the solutions of the 2ODE (\ref{2oder1}) we can write the 1PDE (\ref{eqS}) as a Riccati 1ODE: $ds/dx=s^2+\phi_z\,s-\phi_y$. It is of common knowledge that the transformation $y(x)=-\frac{r'(x)}{f(x)\,r(x)}$ changes the Riccati equation $\,y'(x)=f(x)\,y(x)^2+g(x)\,y(x)+h(x)\,$ into the linear 2ODE $\,r''=(({f'(x)+g(x)\,f(x)})\,r')/{f(x)}-f(x)\,h(x)\,r$. So, with the transformation $s(x)=-w'(x)/w(x)$, the Riccati 1ODE turns (over the solutions of the 2ODE (\ref{2oder1})) into the homogeneous linear 2ODE $\,w''={\phi_z}\,w' + \phi_y\,w$. Therefore, we can use the formal equivalence $D_x\,\sim\,{\frac{d}{dx}}$ to produce a connection between the $S$-functions and the symmetries (written in a particular form) of the 2ODE. Applying the transformation 
\begin{equation}
\label{stonu}
S=-\frac{D_x(\nu)}{\nu}
\end{equation}
into equation (\ref{eqS}). We obtain:
\begin{equation}
\label{2PDEnu}
D_x^2(\nu)=\phi_z\,D_x(\nu)+\phi_y\,\nu.
\end{equation}
The equation (\ref{2PDEnu}) is the symmetry condition for $\,\nu\,$ to be an infinitesimal that defines a symmetry generator in the evolutionary form. So, we have:
\begin{teor}
\label{connectSSy}
Let $\nu$ be a function of $(x,y,z)$ such that $[0,\nu]$ defines a symmetry of the {\em 2ODE (\ref{2oder1})} in the evolutionary form, i.e., $X_e := \nu\,\partial_y\,$ generates a symmetry transformation for  {\em (\ref{2oder1})}. Then the function defined by $S=-D_x(\nu)/\nu$ is a $S$-function associated with the {\em 2ODE (\ref{2oder1})}.
\end{teor}
\noindent
{\it Proof:} See \cite{Noscpc2019}

\begin{cor}
\label{connectSyS}
Let $S$ be a $S$-function associated with the {\em 2ODE (\ref{2oder1})}. Then, the function $\nu$ given by
\begin{equation}
\label{nuS}
\nu \equiv {\rm e}^{\int_x (-S)},
\end{equation}
(where $\, \int_x \,\,\mbox{\rm is the inverse operator of}\,\, D_x\,, \mbox{\rm i.e.}, \, \int_x\,D_x = D_x\,\int_x = \mbox{\boldmath $1$}$) defines a symmetry of the {\em 2ODE (\ref{2oder1})} in the evolutionary form.
\end{cor}
\noindent
{\it Proof:} See \cite{Noscpc2019}

\bigskip

The theorem \ref{connectSSy} and corollary \ref{connectSyS} establish a connection between $S$-functions and Lie symmetries. It is this connection that allowed us to avoid the use of Darboux polynomials in the searching for first integrals for the 2ODE (\ref{2oder1}). The main concept is the {\em associated 1ODE} \footnote{This concept was developed in \cite{Nosjmp2009}, page 222.}, which is an 1ODE that has its general solution defined by one of the first integrals of the 2ODE.

\begin{defin}
\label{1odeass}
Let $I$ be a first integral of the {\em 2ODE (\ref{2oder1})} and let $S(x,y,z)$ be the $S$-function associated with {\em (\ref{2oder1})} through $I$. The first order ordinary differential equation defined by
\begin{equation}
\label{1odeassdefs}
\frac{dz}{dy}= - S(x,y,z),
\end{equation}
where $x$ is taken as a parameter, is called {\bf 1ODE$_{\mathbf{[1]}}$ associated} with {\em (\ref{2oder1})} through $I$.
\end{defin}

\begin{teor}
\label{sol1odeass}
Let $I$ be a first integral of the {\em 2ODE (\ref{2oder1})} and let $S(x,y,z)$ be the $S$-function associated with {\em (\ref{2oder1})} through $I$. Then $I(x,y,z)=C$ is a general solution of the {\em 1ODE$_{[1]}$ associated} with {\em (\ref{2oder1})} through $I$.
\end{teor}
\noindent
{\it Proof:} See \cite{Noscpc2019}

\begin{obs}
\label{solnotinv}
The variable $x$ (the independent variable of the {\em 2ODE (\ref{2oder1})}) is just a parameter in the {\em 1ODE (\ref{1odeassdefs})}, so, since any function of $x$ is an invariant for the operator  $D_a \equiv \partial_y - S\,\partial_z$, i.e., $D_a(x) = 0$, the relation between a general solution $H(x,y,z)=K$ of the {\em 1ODE (\ref{1odeassdefs})} and the first integral $I(x,y,z)$ of the {\em 2ODE (\ref{2oder1})} is $\,I(x,y,z)={\cal F}\left(x,H\right)$, such that $D_x(I) = \partial_x({\cal F}) + (H_x + z\,H_y + \phi\,H_z)\,\partial_H ({\cal F})= 0$.
\end{obs}

In a short way, the $S$-function method consists in determining the $S$-function, solving the associated 1ODE and solving the 1ODE that represents the characteristic system of the partial differential equation for ${\cal F}$. 

\begin{obs}
Some comments:

\begin{itemize}
\item Since the $S$-function is defined as the division of $I_y$ by $I_z$, we can think if the divisions $I_x/I_z$ and $I_x/I_y$ would have any meaning or utility. The answer is yes: We can define the the $S$-functions $\{S_1,S_2,S_3\}$ with similar meanings (see \cite{Noscpc2019}). 
\item The main advantage of the $S$-function method is that, in general, it is more efficient to determine the $S$-function than computing Darboux polynomials precisely in the cases where they (the DPs) have high degrees.
\item A great advantage of having three $S$-functions is that the difficulty associated with determining each one of them can be very different, i.e., it can be much easier to calculate one of them than the others.
\item for more details (and examples) of the application of the method we recommend that the reader see the reference \cite{Noscpc2019}.
\end{itemize}
\end{obs}

\medskip

\section{A method to determine Liouvillian first integrals of LLS equations with an elementary function}
\label{nms1odewf}

In this section we will present a method to find Liouvillian first integrals of LLS equations with an elementary function\footnote{This version of the method will be limited to non algebraic elementary functions.}. In the first subsection we will show an equivalence between a certain type of LLS equations and polynomial vector fields in three variables. In the second, we will show how to use a technique analogous to the one developed in \cite{Noscpc2019} to find a Liouvillian first integral of a polynomial vector field in three variables. In the third subsection, we propose a possible method.  In the fourth subsection, we show how to construct a improvement of the method and, based on it, we propose a second algorithm.

\subsection{LLS equations with functions $\times$ polynomial vector fields}
\label{1odewfpvf}

The method we build is applicable to LLS equations which, after applying the transformation (\ref{tlevsmi}), results in an 1ODE that can be written in the form
\begin{equation}
\label{1odelem}
y' = \frac{dy}{dx} = \phi(x,y,\theta) = \frac{M(x,y,\theta)}{N(x,y,\theta)}
\end{equation}
where  $\phi$ is a rational function of  $(x,y,\theta)$, i.e., $\phi \in K = \C(x,y,\theta)$\footnote{$K = \C(x,y,\theta)$ is the differential field of rational functions in the variables $(x,y,\theta)$.}, $\theta$ is an elementary generator over $\C(x,y)$\footnote{This is a formal way of saying that $\theta$ is an elementary function of $(x,y)$ (see \cite{Dav}).} and $M$ and $N$ are coprime polynomials in $(x,y,\theta)$, i.e., $M,\,N \in \C[x,y,\theta]$\footnote{$\C[x,y,\theta]$ is the differential ring of polynomial functions in the variables $(x,y,\theta)$.}. Furthermore, if $I(x,y,\theta) = c$ represents a general solution of 1ODE (\ref{1odelem}) and $I$ is a Liouvillian function of $(x,y,\theta)$ then for the method to work it is sufficient that the partial derivatives of $I$ in relation to $x$, $y$ and $\theta$ can be expressed as
\begin{eqnarray}
I_x &=& R\,P_1 \label{ix}  \\
I_y &=& R\,P_2 \label{iy}  \\
I_{\theta} &=& R\,P_3 \label{ithe}
\end{eqnarray}
where $R$ is an elementary function of $(x,y,\theta)$, $P_1,\,P_2,\,P_3 \in \C[x,y,\theta]$ (i.e., $P_1,\,P_2$ and $P_3$ are polynomials in $(x,y,\theta)$ ) and the derivatives of $\theta$, $\frac{\partial \theta}{\partial x},\,\frac{\partial \theta}{\partial y} \in K=\C(x,y,\theta)$.

\medskip

The main idea is to perform the transformation $z=\theta(x,y)$ that allows the association of 1ODE (\ref{1odelem}) with a polynomial vector field that has $I(x,y,z)$ as a Liouvillian first integral. Let's see how this happens: 

\begin{itemize}
\item In first place, since $D_x(I)=0$ where $D_x \equiv \partial_x + \phi\,\partial_y$, we have that
\begin{equation}
D_x(I) = \partial_x(I) + \phi\,\partial_y(I) =\frac{\partial I}{\partial x} + \frac{\partial I}{\partial \theta}\,\frac{\partial \theta}{\partial x} + \phi\,\left(\frac{\partial I}{\partial y} + \frac{\partial I}{\partial \theta}\frac{\partial \theta}{\partial y}\right).
\end{equation}

\item Since the derivatives $\frac{\partial \theta}{\partial x},\,\frac{\partial \theta}{\partial y}$ are in $K=\C(x,y,\theta)$, we can write $D_x(I)$ as
\begin{eqnarray}
D_x(I) &=& \frac{\partial I}{\partial x} +  \phi\,\frac{\partial I}{\partial y} + \left(\frac{\partial \theta}{\partial x} + \phi\,\frac{\partial \theta}{\partial y}\right)\,\frac{\partial I}{\partial \theta}=
\nonumber \\ [2mm]
&=& I_x +  \phi\, I_y + \left(\theta_x + \phi\,\theta_y\right)\,I_\theta, \label{eqitem3}
\end{eqnarray}
where $\theta_x \equiv \frac{\partial \theta}{\partial x}$ and $\theta_y \equiv \frac{\partial \theta}{\partial y}$ are rational functions of $(x,y,\theta)$.

\item From the fact that $\phi$, $\theta_x$ and $\theta_y$ are rational functions of $(x,y,\theta)$, we can multiply $D_x(I)$ by the ${\rm lcm}$ of the denominators of $\phi$, $\theta_x$ and $\theta_y\,\phi$ (which we will denote by $l_d$), obtaining 
\begin{equation}
l_d\,D_x(I) = (f\,\partial_x + g\,\partial_y + h\,\partial_{\theta})\left(I(x,y,\theta)\right)=0,
\end{equation}
where $f,\,g$ and $h$ are polynomials in $(x,y,\theta)$ and
\begin{eqnarray}
\label{deefe}
f &=& l_d,
\\ [2mm]
\label{dege}
g &=& l_d\,\phi,
\\ [2mm]
\label{deaga}
h &=& l_d\,\left(\theta_x + \phi\,\theta_y\right).
\end{eqnarray}

\item If we make the substitution $\theta \rightarrow z$, we have that 
\begin{equation}
\label{dii0}
(f(x,y,z)\,\partial_x + g(x,y,z)\,\partial_y + h(x,y,z)\,\partial_z)\left(I(x,y,z)\right)=0.
\end{equation}

\end{itemize}

\medskip

\noindent
Therefore, following these steps we can associate the polynomial vector field 
\begin{equation}
\label{vfchi}
\chi \equiv f(x,y,z)\,\partial_x + g(x,y,z)\,\partial_y + h(x,y,z)\,\partial_z
\end{equation} 
with the 1ODE (\ref{1odelem}). Moreover, we can associate the first integral $I(x,y,z)$ of $\chi$ (see equation (\ref{dii0})) with the general solution of the 1ODE (\ref{1odelem}) (i.e., $I(x,y,\theta)=c$). So, if we can find a Liouvillian first integral $I(x,y,z)$ for the vector field $\chi$, we will have found the general solution of 1ODE (\ref{1odelem}): we only have to apply the substitution $z \rightarrow \theta$ on $I(x,y,z)=c$.

\subsection{$S$-function method adapted to polynomial vector fields in three variables}
\label{sfmpvf3v}

In this section we will show how we can adapt the method described in \cite{Noscpc2019} ($S$-function method) to polynomial vector fields in three variables. Consider that $I(x,y,z)$ is a Liouvillian first integral of the vector field $\chi \equiv f\,\partial_x+g\,\partial_y+h\,\partial_z$, such that the derivatives $I_x$, $I_y$ and $I_z$ are given, respectively, by $R\,P_1$, $R\,P_2$ and $R\,P_3$ (see the section \ref{1odewfpvf}). Now, supose that $I(x,y,z)$ is a Liouvillian first integral of a (hypothetical) rational 2ODE such that $x$ is the independent variable, $y$ is the dependent variable and $z$ is $dy/dx$. This hypothetical 2ODE would be given by 
\begin{equation}
\label{hyp2ode}
z' = \Phi(x,y,z) = - \frac{I_x+z\,I_y}{I_z}= - \frac{P_1+z\,P_2}{P_3} = \frac{M_0(x,y,z)}{N_0(x,y,z)},
\end{equation}
where $M_0$ and $N_0$ are coprime polynomials in $(x,y,z)$.

As we saw in section \ref{sfm}, the $S$-functions $S_1, \, S_2 \, {\rm and} \, S_3$ associated with this 2ODE through the first integral $I$ are given by the following rational functions
\begin{equation}
\label{hypsf}
S_1 = \frac{I_y}{I_z}= \frac{P_2}{P_3}, \,\,\,\,S_2 = \frac{I_x}{I_z}= \frac{P_1}{P_3}, \,\,\,\,S_3 = \frac{I_x}{I_y}= \frac{P_1}{P_2},
\end{equation}
where $S_1, \, S_2 \, {\rm and} \, S_3$ obey the equations
\begin{eqnarray}
\label{s1eq}
 D_x(S_1) &=& {S_1}^2+ \Phi_z\,S_1 - \Phi_y,
 \\ [2mm]
\label{s2eq}
 D_x(S_2) &=& -\frac{1}{z}{S_2}^2+\left(\Phi_z-\frac{\Phi}{z}\right)\,S_2-\Phi_x,
 \\ [2mm]
\label{s3eq}
 D_x(S_3) &=& -\frac{\Phi_y}{\Phi}{S_3}^2+\frac{\Phi_x-z\,\Phi\-y}{\Phi}\,S_3+z\,\Phi_x,
\end{eqnarray}
where $D_x \equiv \partial_x+z\,\partial_y+\Phi\,\partial_z$.
\noindent
Since we do not have the function $\Phi$, we cannot use the equations (\ref{s1eq}), (\ref{s2eq}) and (\ref{s3eq}) to find the $S$-functions $S_1, \, S_2 \, {\rm and} \, S_3$ associated with the hypothetical 2ODE (\ref{hyp2ode}). So, in order to use the technique developed in \cite{Noscpc2019}, we will need some relationship between the (hypothetical) 2ODE and the vector field (\ref{vfchi}). We will establish this relationship using the fact that the 2ODE (\ref{hyp2ode}) and the vector field $\chi$ have in common the first integral $I$. This fact allows us to establish a result that will give us a way to apply the method used in \cite{Noscpc2019} to our problem. Let's begin with the following definition:

\begin{defin}
\label{2ode0}
Let $f\,\partial_x+g\,\partial_y+h\,\partial_z$ be a polynomial vector field in three variables presenting a first integral $I(x,y,z)$ and let $z' = \Phi(x,y,z)$ be a rational 2ODE such that it also has $I(x,y,z)$ as a first integral. We say that the 2ODE and the vector field are {\bf associated through the first integral}~$I$. 
\end{defin}

\begin{obs}
\label{six2odes}
Note that, since the status of the variables $(x,y,z)$ in the vector field $f\,\partial_x+g\,\partial_y+h\,\partial_z$ are the same, there are, in principle, six distinct 2ODEs associated with it:
\begin{eqnarray}
\label{Phi1}
\Phi_1 &=& - \frac{I_x+z\,I_y}{I_z}, \\ [2mm]
\label{Phi2}
\Phi_2 &=& - \frac{I_y+z\,I_x}{I_z}, \\ [2mm]
\label{Phi3}
\Phi_3 &=& - \frac{I_x+y\,I_z}{I_y}, \\ [2mm]
\label{Phi4}
\Phi_4 &=& - \frac{I_z+y\,I_x}{I_y}, \\ [2mm]
\label{Phi5}
\Phi_5 &=& - \frac{I_y+x\,I_z}{I_x}, \\ [2mm]
\label{Phi6}
\Phi_6 &=& - \frac{I_z+x\,I_y}{I_x}. 
\end{eqnarray}
This will be useful in the future process of finding the first integral $I$ through the $S$-function method (see sections \ref{api} and \ref{tpalg}).
\end{obs}

\begin{teor}
\label{s1vfode}
Let $\chi$, defined by {\em (\ref{vfchi})}, be a polynomial vector field presenting a Liouvillian first integral $I$ and let the rational 2ODE defined by {\em (\ref{hyp2ode})} be associated with the vector field $\chi$ through $I$.  If their derivatives are given by 
\begin{eqnarray}
I_x &=& R\,P_1 \label{ix1},  \\
I_y &=& R\,P_2 \label{iy1},  \\
I_z &=& R\,P_3 \label{iz1},
\end{eqnarray}
where $R$ is a Darboux function of $(x,y,z)$ and $P_1,\,P_2,\,P_3 \in \C[x,y,z]$ then the $S$-function $S_1$ associated with the rational 2ODE {\em (\ref{hyp2ode})} is given by  
\begin{equation}
\label{s12ode1}
S_1 = \frac{M_0\,f - N_0\,h}{N_0\,(g-z\,f)},
\end{equation}
where $g-z\,f \neq 0$.
\end{teor}

\noindent
{\bf Proof of Theorem \ref{s1vfode}:} From the hypotheses we have that

\begin{eqnarray}
\label{dosi}
D_0(I) & = & N_0 \,\partial_x(I) + z\, N_0 \,\partial_y(I)\, + M_0 \,\partial_z(I) = 0,  \\ [2mm]
\chi(I) & = & f \,\partial_x(I) + g \,\partial_y(I)\, + h \,\partial_z(I) = 0. \label{dosi2}
\end{eqnarray}
where $M_0$ and $N_0$ are, respectively, the numerator and denominator of the 2ODE (\ref{hyp2ode}). Defining $D_1 \equiv f\,D_0 - N_0\,\chi$, and applying it to $I$ we obtain:
\begin{equation}
\label{d1i}
D_1 (I) = ( f\,D_0 - N_0\,\chi) (I) = (z\,f-g)\,N_0 \,\partial_y(I) + (M_0\,f - N_0\,h)\,\partial_z(I) = 0.
\end{equation}
From eq.(\ref{d1i}) we have (see \cite{Noscpc2019}) that $I_y/I_z=-(M_0\,f - N_0\,h)/(N_0\, (z\,f-g))=S_1. \,\,\,\Box$

\begin{obs}
We note that the function $S_1$ as presented in the equation {\em (\ref{s12ode1})} `defines', in a certain sense, the relationship between the vector field $\chi$ and the 2ODE {\em (\ref{hyp2ode})}. So, we can use it to produce an algorithm (semi) to find the 2ODE ($M_0$ and $N_0$) and, in view of {\em (\ref{s12ode1})}, the $S$-function.
\end{obs}

\bigskip

\subsection{A possible method}
\label{mffi}

In this section, we will use the equation (\ref{s12ode1}) in the process of building an algorithm to find the general solution of 1ODE (\ref{1odelem}). Before describing the steps of a possible algorithm in more detail, let's outline the main steps of the process and discuss its application to some examples in order to materialize (and clarify) the path we are trying to follow to achieve our goal.

\medskip

\noindent
{\bf The Procedure:} (Sketch)
\begin{itemize}
\item If the 1ODE can be placed in the form (\ref{1odelem}) then construct the operator $D$ (\ref{eqitem3}).
\item Make the substitution $\theta \rightarrow z$ on the $D$ operator and find the polynomial operator $\chi$ (\ref{vfchi}).
\item Construct polynomial candidates $M_c$ and $N_c$ (with undetermined coefficients $\{m_i\}$ and $\{n_j\}$) of some chosen degree in the variables $(x,y,z)$.
\item Substitute 

\noindent
$S_1 = \displaystyle{\frac{M_c\,f - N_c\,h}{N_c\,(g-z\,f)}}$

\noindent
in the equation for the $S$-function: 

\noindent
$D_x(S_1) = {S_1}^2+ \displaystyle{\frac{\partial \Phi_c}{\partial z}\,S_1 - \frac{\partial \Phi_c}{\partial y}},$

\noindent 
where $D_x \equiv \partial_x+z\,\partial_y+\Phi_c\,\partial_z$ and $\Phi_c \equiv M_c/N_c$.
\item Solve the polynomial equation for the unknown coefficients $\{m_i\}$ and $\{n_j\}$, obtaining $S_1$.
\item Use the $S$-function method to find de Liouvillian first integral $I$ of the 2ODE $z' = \Phi(x,y,z)$.
\item Make the substitution $z \rightarrow \theta$ on the first integral $I$ and obtain the desired general solution of the 1ODE.
\end{itemize}

\medskip

Before presenting a more formal algorithm let's discuss some examples:

\noindent
{\bf Example 1:} Consider the 1ODE given by
\begin{equation}
\label{exem1}
\frac{dy}{dx} = {\frac {{{\rm e}^{x}}{x}^{3}{y}^{2}+{{\rm e}^{x}}{x}^{2}{y}^{2}+2\,{
{\rm e}^{x}}{x}^{2}y+{{\rm e}^{x}}xy+{{\rm e}^{x}}x+{y}^{2}+{{\rm e}^{
x}}}{{x}^{2}{y}^{2}+{{\rm e}^{x}}{x}^{2}+xy+1}}
\end{equation}
and let's apply the procedure sketched above:

\begin{itemize}

\item Choosing $\theta = {\rm e}^x$, we have that $\theta_x = {\rm e}^x = \theta$ and $\theta_y = 0$. So, the operator $D$ will be:
\begin{equation}
\label{ex1D}
D = \partial_x + {\frac {\theta\,{x}^{3}{y}^{2}+\theta\,{x}^{2}{y}^{2}+2\,\theta\,{x}^{2}y+\theta\,xy+\theta\,x+{y}^{2}+\theta}{{x}^{2}{y}^{2}+\theta\,{x}^{2}+xy+1}}\,\partial_y + \theta\,\partial_{\theta}.
\end{equation}

\item The polynomials $f$, $g$ and $h$ will be:
\begin{eqnarray}
\label{ex1f}
f &=& {x}^{2}{y}^{2}+z{x}^{2}+xy+1,
 \\ [2mm]
\label{ex1g}
 g &=& z{x}^{3}{y}^{2}+z{x}^{2}{y}^{2}+2\,z{x}^{2}y+zxy+zx+{y}^{2}+z,
 \\ [2mm]
\label{ex1h}
 h &=& z \left( {x}^{2}{y}^{2}+z{x}^{2}+xy+1 \right).
\end{eqnarray}

The operator $\chi$ will be:
\begin{equation}
\label{ex1chi}
\chi = f\, \partial_x + g\,\partial_y + h\,\partial_{z}.
\end{equation}

\item For the degrees ${\rm deg}_M=4$ and ${\rm deg}_N=5$ we obtain:
\begin{eqnarray}
\label{ex1M}
M_0 &=& -\left( xz-y \right)  \left( xz+y \right),
 \\ [2mm]
\label{ex1N}
 N_0 &=& - x \left( xy+1 \right) ^{2},
 \\ [2mm]
\label{ex1S1}
 S_1 &=& - \frac{{x}^{2}{y}^{2}+{x}^{2}z+xy+1}{x \left( xy+1 \right) ^{2}}
\end{eqnarray}

\item Using $S_1$ in the $S$-function method, we have the following first integral:
\begin{equation}
\label{ex1I}
I = {\rm e}^{\frac{1}{xy+1}} \left( zx-y \right)
\end{equation}

\item Making the substitution $z={\rm e}^x$ we finally arrive at the solution:
\begin{equation}
\label{ex1sol}
{\rm e}^{{\frac{1}{xy+1}}} \left( x\,{\rm e}^x -y \right) = C
\end{equation}

\end{itemize}

\noindent
Now let's take a closer look at the derivatives of the first integral $I$ (of the vector field $\chi$):
\begin{eqnarray}
I_x &=& \frac {{\rm e}^{\frac{1}{xy+1}}}{(xy+1)^2}\,\left( z{x}^{2}{y}^{2}+
zxy+{y}^{2}+z \right), 
\label{ixex1} \\ [2mm]
I_y &=& - \frac {{\rm e}^{\frac{1}{xy+1}}}{(xy+1)^2} \, \left( {x}^{2}{y}^{2}+
z{x}^{2}+xy+1 \right), 
\label{iyex1} \\ [2mm]
I_z &=& \frac {{\rm e}^{\frac{1}{xy+1}}}{(xy+1)^2} \, x (xy+1)^2. 
\label{izex1}
\end{eqnarray}
From them we can determine (directly) the integrating factor $R$, the 2ODE associated ($\Phi$), the polynomials $P_1$, $P_2$ and $P_3$ and the $S$-functions $S_1$, $S_2$ and $S_3$: 
\begin{eqnarray}
\Phi &=& {\frac {{x}^{2}{z}^{2}-{y}^{2}}{ \left( xy+1 \right) ^{2}x}},
\label{Phiex1} \\ [2mm]
R &=& \frac {{\rm e}^{\frac{1}{xy+1}}}{(xy+1)^2},
\label{rex1} \\ [2mm]
P_1 &=& {x}^{2}{y}^{2}z+xyz+{y}^{2}+z,
\label{p1ex1} \\ [2mm]
P_2 &=& - \left( {x}^{2}{y}^{2}+{x}^{2}z+xy+1 \right), 
\label{p2ex1} \\ [2mm]
P_3 &=&  x (xy+1)^2, 
\label{p3ex1} \\ [2mm]  
S_1 &=& -{\frac {{x}^{2}{y}^{2}+z{x}^{2}+xy+1}{ \left( xy+1 \right) ^{2}x}},
\label{s1ex1} \\ [2mm]
S_2 &=&{\frac {z{x}^{2}{y}^{2}+zxy+{y}^{2}+z}{ \left( xy+1 \right)^{2}x}}, 
\label{s2ex1} \\ [2mm]
S_3 &=&  -{\frac {z{x}^{2}{y}^{2}+zxy+{y}^{2}+z}{{x}^{2}{y}^{2}+z{x}^{2}+xy+1}}.
\label{s3ex1} 
\end{eqnarray}
The components $(f,g,h)$ of the vector field $\chi$ are (see above)
\begin{eqnarray}
f &=& {x}^{2}{y}^{2}+z{x}^{2}+xy+1,
 \\ [2mm]
g &=& z{x}^{3}{y}^{2}+z{x}^{2}{y}^{2}+2\,z{x}^{2}y+zxy+zx+{y}^{2}+z,
 \\ [2mm]
h &=& z \left( {x}^{2}{y}^{2}+z{x}^{2}+xy+1 \right).
\end{eqnarray}

\begin{obs}
\label{ex1obs}
So, comparing $(f,g,h)$ with the polynomials $P_1$, $P_2$ and $P_3$ and with the numerators and denominators of the $S$-functions, we can see that:
\begin{equation}
N_0 = P_3,\,\,\,\, f = - P_2,\,\,\,\, h = - z\,P_2.
\end{equation}
The first of these equations may not surprise, since $\Phi = -(P_1+z\,P_2)/P_3$. However, the following two equations are (at first glance) intriguing. 
\end{obs}

\medskip

\noindent
{\bf Example 2:} Consider now the 1ODE given by
\begin{equation}
\label{exem2}
\frac{dy}{dx} ={\frac {{-{\rm e}^{y}} \left( {{\rm e}^{2y}}{x}^{2}y-2\,{{\rm e}^{y}}x{y}^{2}-x{{\rm e}^{y}}y+{y}^{3}-1 \right) }
{ {{\rm e}^{3y}}{x}^{3}y\!+\! {{\rm e}^{3y}}{x}^{3}\!-2{{\rm e}^{2y}}{x}^{2}{y}^{2}\!-3{{\rm e}^{2y}}{x}^{2}y\!+{{\rm e}^{y}}x{y}^{3}\!
+{{\rm e}^{y}}x{y}^{2}\!+x{{\rm e}^{y}}y\!-x{{\rm e}^{y}}\!+\!1}}.
\end{equation}
Applying the whole procedure again we obtain 

\begin{eqnarray}
\label{ex2f}
f &=&\! {x}^{3}y{z}^{3}\!+\!{z}^{3}{x}^{3}\!-\!2{x}^{2}{y}^{2}{z}^{2}\!-\!3{z}^{2}{x}^{2}y\!+\!x{y}^{3}z\!+\!zx{y}^{2}\!+\!xzy\!-\!zx\!+\!1,
 \\ [2mm]
\label{ex2g}
 g &=& -z \left( {z}^{2}{x}^{2}y-2\,zx{y}^{2}-xzy+{y}^{3}-1 \right),
 \\ [2mm]
\label{ex2h}
 h &=& -{z}^{2} \left( {z}^{2}{x}^{2}y-2\,zx{y}^{2}-xzy+{y}^{3}-1 \right).
\end{eqnarray}

\noindent
The function $S_1$ is
\begin{equation}
\label{ex2S1}
 S_1 = {\frac {{z}^{3}{x}^{3}-2\,{z}^{2}{x}^{2}y+zx{y}^{2}+xzy+1}{x \left( {z
}^{2}{x}^{2}y-2\,zx{y}^{2}-xzy+{y}^{3}-1 \right) }}.
\end{equation}

\noindent
From $S_1$ we can determine the first integral (for the vector field $\chi$) and (making the substitution $z={\rm e}^y$) the general solution of the 1ODE (\ref{exem2}):
\begin{equation}
\label{ex2Isol}
I = {\rm e}^{\frac{1}{zx-y}} \left( xzy+1 \right), \,\,\,\, {\rm e}^{\frac{1}{{\rm e}^{y}x-y}} \left( x{{\rm e}^{y}}y+1 \right)=C.
\end{equation}

\medskip

\noindent
From the derivatives of the first integral $I$ we can get $\Phi$, the polynomials $P_1$, $P_2$ and $P_3$ and the $S$-functions $S_1$, $S_2$ and $S_3$: 
\begin{eqnarray}
\Phi &=& {\frac { \left( zx-y \right) ^{2} \left( zx+y \right) z}{x \left( {z}^
{2}{x}^{2}y-2\,zx{y}^{2}-xzy+{y}^{3}-1 \right) }},
\label{Phiex2} \\ [2mm]
P_1 &=& z \left( {z}^{2}{x}^{2}y-2\,zx{y}^{2}-xzy+{y}^{3}-1 \right),
\label{p1ex1} \\ [2mm]
P_2 &=& {z}^{3}{x}^{3}-2\,{z}^{2}{x}^{2}y+zx{y}^{2}+xzy+1, 
\label{p2ex1} \\ [2mm]
P_3 &=&  x \left( {z}^{2}{x}^{2}y-2\,zx{y}^{2}-xzy+{y}^{3}-1 \right), 
\label{p3ex1} \\ [2mm]  
S_1 &=& {\frac {{z}^{3}{x}^{3}-2\,{z}^{2}{x}^{2}y+zx{y}^{2}+xzy+1}{x \left( {z
}^{2}{x}^{2}y-2\,zx{y}^{2}-xzy+{y}^{3}-1 \right) }},
\label{s1ex1} \\ [2mm]
S_2 &=&{\frac {z}{x}}, 
\label{s2ex1} \\ [2mm]
S_3 &=&  {\frac {z \left( {z}^{2}{x}^{2}y-2\,zx{y}^{2}-xzy+{y}^{3}-1 \right) }{
{z}^{3}{x}^{3}-2\,{z}^{2}{x}^{2}y+zx{y}^{2}+xzy+1}}.
\label{s3ex1} 
\end{eqnarray}

\medskip

\begin{obs}
\label{ex2obs}
Again, if we compare $(f,g,h)$ with the polynomials $P_1$, $P_2$ and $P_3$ and with the $S$-functions, we see that:
\begin{equation}
N_0 = P_3,\,\,\,\, g = - P_1,\,\,\,\, h = - z\,P_1.
\end{equation}
Besides that, as the roles of $x$ and $z$ can be switched in the first integral $I$ (i.e., the transformation $x \rightarrow z, \, z\rightarrow x$ is a symmetry transformation for $I$) the $S$-function $S_2$ has a very simple format.
\end{obs}

\medskip

{\bf Example 3:} Let's see now an example where we have a function of $x$ and $y$. Consider the 1ODE
\begin{equation}
\label{exem3}
\frac{dy}{dx} =\frac{xy\ln  \left( {\frac {x}{y}} \right) + \left( x{y}^{5}-2\,{x}^{2}{y}^{
3}+{x}^{3}y+{y}^{4}+{x}^{2}y-2\,x{y}^{2}+{x}^{2} \right) y}
{2\,x{y}^{2}\ln  \left( {\frac {x}{y}} \right) +x \left( -x{y}^{5}+2\,{
x}^{2}{y}^{3}-{x}^{3}y+2\,x{y}^{3}+{y}^{4}-2\,x{y}^{2}+{x}^{2} \right)}.
\end{equation}

\noindent
So
\begin{eqnarray}
\label{ex3f}
\!\!\!\!f\! &=&\!\! x \left( x{y}^{5}-2\,{x}^{2}{y}^{3}+{x}^{3}y-2\,x{y}^{3}-{y}^{4}+2\,x{y}^{2}-2\,z{y}^{2}-{x}^{2} \right),
 \\ [2mm]
\label{ex3g}
\!\!\!\! g\! &=&\! - \left( x{y}^{5}-2\,{x}^{2}{y}^{3}+{x}^{3}y+{y}^{4}+{x}^{2}y-2\,x{y}^{2}+{x}^{2}+zx \right) y,
 \\ [2mm]
\label{ex3h}
\!\!\!\! h\! &=&\! 2\,x{y}^{5}-4\,{x}^{2}{y}^{3}+2\,{x}^{3}y-2\,x{y}^{3}+{x}^{2}y-2\,z{y}^{2}+zx.
\end{eqnarray}

\noindent
The function $S_1$ is
\begin{equation}
\label{ex3S1}
 S_1 = {\frac {x{y}^{4}-2\,{x}^{2}{y}^{2}+{x}^{3}-2\,x{y}^{2}-2\,yz}{{y}^{4}-
2\,x{y}^{2}+{x}^{2}}}.
\end{equation}

\noindent
The first integral and the general solution of the 1ODE are
\begin{equation}
\label{ex3Isol}
I = {\rm e}^{\frac{1}{y^2-x}} \left( xy+z \right), \,\,\,\, {\rm e}^{\frac{1}{y^2-x}} \left( xy+\ln  \left( {\frac {x}{y}} \right) \right)=C.
\end{equation}

\noindent
The polynomials $P_1$, $P_2$, $P_3$ and $\Phi$ are: 
\begin{eqnarray}
\!\!\!P_1\!\! &=&\!\! {y}^{5}-2\,x{y}^{3}+{x}^{2}y+xy+z,
\label{p1ex3} \\ [2mm]
\!\!\!P_2\!\! &=&\!\! x{y}^{4}-2\,{x}^{2}{y}^{2}+{x}^{3}-2\,x{y}^{2}-2\,yz, 
\label{p2ex3} \\ [2mm]
\!\!\!P_3\!\! &=&\!\! \left( -{y}^{2}+x \right) ^{2}, 
\label{p3ex3} \\ [2mm]
\!\!\!\Phi\!\! &=&\!\! {\frac {x{y}^{4}z\!-\!2{x}^{2}{y}^{2}z\!+\!{y}^{5}\!+\!{x}^{3}z\!-\!2x{y}^{3}\!-\!2x
{y}^{2}z\!+\!{x}^{2}y\!-\!2y{z}^{2}\!+\!xy\!+\!z}{ \left( -{y}^{2}+x \right) ^{2}}}.
\label{Phiex3}
\end{eqnarray}

\begin{obs}
\label{ex3obs}
Some comments:
\begin{enumerate} 
\item This time, looking at the polynomials $P_1$, $P_2$ and $P_3$ in the example 3, we see that none of them is a member of $\{c\,f,c\,g,c\,h\}$ where $c$ is a constant.

\item A very noticeable difference between examples 1 and 2 and example 3 is that in examples 1 and 2 the elementary function present on the 1ODE was a function of only one variable ($x$ or $y$) whereas in example 3 the function $\theta$ was a function of the two variables $(x,y)$, i.e., $\theta=\theta(x,y)$.

\item In general it is much simpler (computationally speaking) to calculate the $S$-function if we already know its numerator or denominator. In these cases the algorithm is much more efficient (i.e., faster and with less memory consumption). 

\item What happened in examples 1 and 2 was not a fluke. We can show that, in almost all cases, if $\theta$ is an elementary function of only one variable, then one of the coefficients of the vector field $\chi$ will divide one of the polynomials $(P_1,P_2,P_3)$. 

\item The good news is that, if $\theta$ is an elementary function of a rational function, in a lot of cases we can perform a variable transformation that takes $\theta(x,y)$ into $\overline{\theta}(x)$ (or $\overline{\theta}(y)$). 
\end{enumerate} 
\end{obs}

\subsection{A possible improvement}
\label{api}

What was stated in observation 4 (remark \ref{ex3obs}) can be demonstrated and (in conjunction with the fact presented in observation 5) used to greatly improving the efficiency of the algorithm in a wide variety of cases\footnote{However, although the improvement we are going to present is quite general, the `problem solution' can, in certain cases, introduce complexity to the process, almost canceling the efficiency gain obtained.}. Before presenting the alternative method, let's define (a little more) the type of 1ODE that these improvement can handle.

\begin{defin}
\label{lsset}
Consider that an 1ODE as described by equation {\em (\ref{1odelem})} has the following characteristics:
\begin{enumerate}
\item  $I(x,y,\theta(x,y)) = c$ ($c$ constant) represents a general solution of the 1ODE {\em (\ref{1odelem})} and $I$ is a Liouvillian function of $(x,y,\theta)$.
\item The function $\theta(x,y)$ is in an elementary extension $E$ of the differential field (of characteristic zero) $\C(x,y)$ such that $E=\C(x,y,\exp(r))$ or $E=\C(x,y,\ln(r))$, where $r \in \C(x,y)$.
\item The 1ODE is associated with a vector field $\chi \equiv f\,\partial_x+g\,\partial_y+h\,\partial_z$ such that $f,g,h \in \C[x,y,z]$ and
\begin{eqnarray}
f &=& l_d |_{\theta \rightarrow z}, \nonumber 
\\ [2mm]
g &=& l_d\,\phi |_{\theta \rightarrow z}, \nonumber 
\\ [2mm]
h &=& l_d\,\left(\theta_x + \phi\,\theta_y\right) |_{\theta \rightarrow z}. \nonumber
\end{eqnarray}
\item The derivatives of the first integral $I$ (of the vector field $\chi$) are of the form: $\partial_k (I) = \exp(A/B)\,\prod_i p_i^{n_i}\,P_k,\,k \in \{1,2,3\}$, where $A,B,p_i,P_k$ are polynomials in $\C[x,y,z]$ and $n_i$ are integers. 
\item Let $T$ be the transformation $T=\{u=r(x,y),v=y\}$ (or $T = \{u = x, v = r(x,y) \}$), where $r \in \C(x,y)$ is the argument of $\theta$. If $T^{-1} = \{x=r_1(u,v),y=v\}$ (or $T^{-1} = \{ x = u, y=r_1(u,v) \}$) denote its inverse, we are assuming that $r_1$ is a rational function of $(u,v)$, i.e., $r_1 \in \C(u,v)$.
\end{enumerate}
We will denote the set of 1ODEs that have these characteristics as $L_S$.
\end{defin}

Now, we are going to present a result that, if applied to an 1ODE $\in L_S$, can be used to construct an improvement of the method sketched above.

\begin{teor}
\label{theoAimp}
Let $\chi \equiv f\,\partial_x+g\,\partial_y+h\,\partial_z$ be the polynomial vector field associated with the 1ODE {\em (\ref{1odelem})} ($\phi = M/N,\,M,N \in   \C[x,y,\theta]$) and let the first integral $I$ (of the vector field $\chi$) be such that its derivatives are of the form: $\partial_k (I) = \exp(A/B)\,\prod_i p_i^{n_i}\,P_k,\,k \in \{1,2,3\}$, where $A,B,p_i,P_k$ are polynomials in $\C[x,y,z]$ and $n_i$ are rational numbers. We can state that:

\noindent
i) If $\,\theta = {\rm e}^x\,$ then $\,f | P_2,\,$ $h=z\,f\,$  and $\,g |(P_1 + z\,P_3)$.

\noindent
ii) If $\,\theta = \ln(x)\,$ then $\,h | P_2,\,$ $f=x\,h\,$ and $\,g | (x\,P_1+P_3)$.

\noindent
iii) If $\,\theta = {\rm e}^y\,$ then $\,g | P_1,\,$ $h=z\,g\,$ and $\,f |(P_2 + z\,P_3)$.

\noindent
iv) If $\,\theta = \ln(y)\,$ then $\,h | P_1,\,$ $g=y\,h\,$ and $\,f | (y\,P_2+P_3)$.
\end{teor}

\noindent
{\it Proof.} We have seen that we can write the operator $D \equiv \partial_x + \phi\,\partial_y$ -- associated with the 1ODE (\ref{1odelem}) -- in the form
$
D ={\partial_x} +  \phi\,{\partial_y} + \left(\theta_x + \phi\,\theta_y\right)\,{\partial_\theta},
$
where $\theta_x$ and $\theta_y$ are rational functions of $(x,y,\theta)$.
In order to obtain the vector field $\chi$ we have to multiply $D$ by the ${\rm lcm}$ of the denominators of $\phi$, $\theta_x$ and $\theta_y\,\phi$ and make the substitution $\theta \rightarrow z$, obtaining  
$
\chi = f\,\partial_x + g\,\partial_y + h\,\partial_z, \nonumber 
$
where $f,\,g$ and $h$ are polynomials given by
\begin{eqnarray}
f &=& l_d |_{\theta \rightarrow z}, \nonumber 
\\ [2mm]
g &=& l_d\,\phi |_{\theta \rightarrow z}, \nonumber 
\\ [2mm]
h &=& l_d\,\left(\theta_x + \phi\,\theta_y\right) |_{\theta \rightarrow z}. \nonumber
\end{eqnarray}

\noindent
Now we can prove the statements: 

\medskip

\noindent
$i$): If $\theta={\rm e}^x$ then $\theta_x = z$ and $\theta_y = 0$ and so
\begin{equation}
\label{proo1}
f=N(x,y,z),\,\,g=M(x,y,z), \,\,h=z\,N(x,y,z).
\end{equation}
We have that $\chi(I)\! =\! f\,I_x+g\,I_y+h\,I_z= fRP_1+gRP_2+hRP_3 \,{\rm (by\, hypothesis)} = R\, (f\,P_1+g\,P_2+h\,P_3) = 0$. Substituting (\ref{proo1}) and noting that $R$ is not null we can write $N\,P_1+M\,P_2+z\,N\,P_3=0$ implying that $N\,(P_1+z\,P_3)=-M\,P_2$. Therefore, we can write 
\begin{equation}
\label{proo2}
P_1+z\,P_3=- \frac{M\,P_2}{N} \,\, {\rm and} \,\, P_2=- \frac{N\,(P_1+z\,P_3)}{M}.
\end{equation}
Since $M$ and $N$ are coprime, we have that $\,N|P_2\,\,\Rightarrow\,\,f | P_2\,$ and $\,M|(P_1 + z\,P_3)\,\,\Rightarrow\,\,g |(P_1 + z\,P_3)$.

\medskip

\noindent
$ii$): If $\theta=\ln(x)\,\,\Rightarrow\,\,\theta_x = 1/x$ and $\theta_y = 0$. We have two cases:
\begin{equation}
\label{proo3}
f=x\,N(x,y,z),\,\,g=x\,M(x,y,z), \,\,h=N(x,y,z),
\end{equation}
if $x$ is not a factor of $N$ or 
\begin{equation}
\label{proo4}
f=x^n\,P_N(x,y,z),\,\,g=M(x,y,z), \,\,h=x^{n-1}P_N(x,y,z),
\end{equation}
where $P_N=N/x^n$ is a polynomial, $n$ is a positive integer and $x$ and $P_N$ are coprime. 

\noindent
First case ($x$ and $N$ are coprime): $N\,(x\,P_1+P_3)=-x\,M\,P_2$. Therefore, $\,N|P_2\,\,\Rightarrow\,\,h | P_2\,$ and $\,x\,M|(x\,P_1+P_3)\,\,\Rightarrow\,\,g |(x\,P_1+P_3)$. 

\noindent
Second case ($x|N$): $x^{n-1}P_N\,(x\,P_1+P_3)=-\,M\,P_2$. Since $M$ and $N$ are coprime, we have that $x^{n-1}P_N\,(=h)$ and $M$ are coprime. Therefore, $h | P_2\,$ and $g |(x\,P_1+P_3)$.

\medskip

\noindent
$iii$): If $\theta={\rm e}^y$ then $\theta_x = 0$ and $\theta_y = z$ and so
\begin{equation}
\label{proo5}
f=N(x,y,z),\,\,g=M(x,y,z), \,\,h=z\,M(x,y,z).
\end{equation}
We can write $N\,P_1+M\,P_2+z\,M\,P_3=0$ implying that $N\,P_1=-M\,(P_2+z\,P_3)$. Therefore, we can write 
\begin{equation}
\label{proo6}
P_2+z\,P_3=- \frac{N\,P_1}{M} \,\, {\rm and} \,\, P_1=- \frac{M\,(P_2+z\,P_3)}{N}.
\end{equation}
Since $M$ and $N$ are coprime, we have that $\,N|(P_2+z\,P_3)\,\,\Rightarrow\,\,f | (P_2+z\,P_3)\,$ and $\,M|P_1 \,\,\Rightarrow\,\,g | P_1$.

\medskip

\noindent
$iv$): If $\theta=\ln(y)\,\,\Rightarrow\,\,\theta_x = 0$ and $\theta_y = 1/y$. We have two cases:
\begin{equation}
\label{proo7}
f=y\,N(x,y,z),\,\,g=y\,M(x,y,z), \,\,h=M(x,y,z),
\end{equation}
if $y$ is not a factor of $M$ or 
\begin{equation}
\label{proo8}
f=N(x,y,z),\,\,g=y^{n}P_M(x,y,z), \,\,h=y^{n-1}P_M(x,y,z),
\end{equation}
where $P_M=M/y^n$ is a polynomial, $n$ is a positive integer and $y$ and $P_M$ are coprime. 

\noindent
First case ($y$ and $M$ are coprime): $N\,y\,P_1=-M\,(y\,P_2+P_3)$. Therefore, $\,M|P_1\,\,\Rightarrow\,\,h | P_1\,$ and $\,y\,N|(y\,P_2+P_3)\,\,\Rightarrow\,\,f |(y\,P_2+P_3)$. 

\noindent
Second case ($y|M$): $y^{n-1}P_M\,(y\,P_2+P_3)=-\,N\,P_1$. Since $M$ and $N$ are coprime, we have that $y^{n-1}P_M\,(=h)$ and $N$ are coprime. Therefore, $h | P_1\,$ and $f |(y\,P_2+P_3)$. $\,\,\,\Box$

\medskip

\begin{cor}
\label{coroAimp}
If the 1ODE {\em (\ref{1odelem})} $\in L_S$, then there is a rational transformation $T$ (with rational inverse $T^{-1}$) such that the transformed 1ODE$_{[tr]}$ has an associated vector field $\chi_{[tr]}$ in which $f_i|P_j$ for some $\,f_i \in \{f,g,h\}\,$ and for some $\,P_j \in \{P_1,P_2\}$. 
\end{cor}

\noindent
{\it Proof.} For any 1ODE $\in L_S$ we can perform a variable transformation that takes $\theta(x,y)$ into $\exp(x)$, $\ln(x)$, $\exp(y)$ or $\ln(y)$. Thus, any 1ODE $\in L_S$ can be transformed into an 1ODE that fulfills the premises of the theorem \ref{theoAimp}. Therefore, the conclusion is a direct consequence of the theorem \ref{theoAimp}. $\,\,\,\Box$

\medskip

In this point we can use the compatibility conditions for the first integral $I$ of the vector field $\chi$ to write the 1PDE for $S$ in terms of $f,\, g$ and $h$.

\begin{teor}
\label{theoSsis}
Let $\chi \equiv f\,\partial_x+g\,\partial_y+h\,\partial_z$ be the polynomial vector field associated with the 1ODE {\em (\ref{1odelem})} and let $I$ be a first integral it. Then the $S$-function associated with $\chi$ through $I$ obeys the 1PDE given by
\begin{equation}
\label{eqSsys}
\chi (S)=S^2\left(\frac{fg_z-gf_z}{f}\right)+S\left(\frac{gf_y-fg_y+fh_z-hf_z}{f}\right)-\left(\frac{fh_y-hf_y}{f}\right).
\end{equation}
\end{teor}

\noindent
{\it Proof.} See the proof of theorem 1.1 and corollary 1.1 of \cite{Noscpc2019}. $\,\,\,\Box$

\medskip

\noindent
Since $S=I_y/I_z=P_2/P_3$ and the status of the variables $(x,y,z)$ in the vector field $f\,\partial_x+g\,\partial_y+h\,\partial_z$ is the same, we can carry out a procedure consisting of the following steps:

\medskip

{\bf Procedure} (sketch):
\begin{itemize}

\item In first place we can use the transformation $T_1$ such that its inverse is given by ${T_1}^{-1}=\{y_1=x,x_1=r(x,y)\}$ (where $r(x,y)$ is the argument of $\theta$) to obtain an 1ODE$_{[1]}$ such that $\theta_{[1]} = \exp(x_1)$ or $\theta_{[1]} = \ln(x_1)$.

\item If $\theta_{[1]} = \exp(x_1)$ we can use the transformation $T_2=\{x_1=\ln(x_2),y_1=y_2\}$ and obtain an 1ODE$_{[2]}$ such that $\theta_{[2]} = \ln(x_2)$.

\item For the 1ODE$_{[2]}$, the vector field $\chi_2$ associated is given by $f_2\,\partial_{x_2}+g_2\,\partial_{y_2}+h_2\,\partial_{z_2}$, where (To make the notation less heavy we omitted index 2 in the following equations):
\begin{equation}
f=x\,N(x,y,z),\,\,g=x\,M(x,y,z), \,\,h=N(x,y,z),
\end{equation}
if $x$ is not a factor of $N$ or 
\begin{equation}
f=x^n\,P_N(x,y,z),\,\,g=M(x,y,z), \,\,h=x^{n-1}P_N(x,y,z),
\end{equation}
where $P_N=N/x^n$ is a polynomial, $n$ is a positive integer and $x$ and $P_N$ are coprime (see case ($ii$) of theorem \ref{theoAimp}). 

\item Since the status of the variables $(x_2,y_2,z_2)$ in the vector field $\chi_2$ is the same, we can make the transformation $T_3=\{x_2=y_3,y_2=z_3,z_2=x_3\}$ that will make $P_2$ take the place of $P_3$ (which is the denominator of the transformed $S$-function). The transformed vector field $\chi_3$ will be given by $f_3\,\partial_{x_3}+g_3\,\partial_{y_3}+h_3\,\partial_{z_3}$, where
\begin{eqnarray}
f_3 &=& T_3(h_2), \nonumber 
\\ [2mm]
g_3 &=& T_3(f_2), \nonumber 
\\ [2mm]
h_3 &=& T_3(g_2). \nonumber
\end{eqnarray}

\item At this point we can use the equation (\ref{eqSsys}) to determine the numerator of the transformed $S$-function $S_{1\,[3]}$. 

\item We can use $S_{1\,[3]}$ to find a first integral $I_{[3]}$ for the vector field $\chi_{[3]}$. 

\item Applying the inverse trasformation ${T_3}^{-1}$ to $I_{[3]}$ we obtain the first integral $I_{[2]}$ for $\chi_{[2]}$.  

\item From $I_{[2]}$ we can write the general solution of 1ODE$_{[2]}$ (using $z\,\rightarrow\,\ln(x)$). 

\item If we apply the inverse trasformation ${T_2}^{-1}$ followed by ${T_1}^{-1}$ to $I_{[2]}=c$, we will obtain $I = {T_1}^{-1}({T_2}^{-1}(I_{[2]}))=c$, which represents the general solution of the original 1ODE. 
\end{itemize}

Let's make these steps clearer with an example:

\bigskip

\noindent
{\bf Example 4:} Consider the 1ODE given by
\begin{equation}
\label{exem4}
\frac{dy}{dx} = {\frac {2\, \left( {{\rm e}^{xy}} \right) ^{2}{x}^{3}{y}^{2}-{{\rm e}^
{xy}}{x}^{3}{y}^{2}+{{\rm e}^{xy}}x{y}^{3}-5\,{{\rm e}^{xy}}{x}^{2}y+{
{\rm e}^{xy}}{y}^{2}+2\,x}{ \left( {{\rm e}^{xy}} \right) ^{2}{x}^{2}{
y}^{2}+{{\rm e}^{xy}}{x}^{4}y-{{\rm e}^{xy}}{x}^{2}{y}^{2}+{{\rm e}^{x
y}}{x}^{3}-3\,xy{{\rm e}^{xy}}+1}}.
\end{equation}

\begin{itemize}

\item In first place, let's make the transformation $T_1 = \{x=x/y,\,y=y\}$ leading to:
\begin{equation}
\label{1ode1ex4}
\!\!\!\!\frac{dy}{dx} \!=\! {\frac {y \left( 2\,{{\rm e}^{2\,x}}{x}^{3}+{{\rm e}^{x}}x{y}^{3}-{
{\rm e}^{x}}{x}^{3}+{{\rm e}^{x}}{y}^{3}-5\,{{\rm e}^{x}}{x}^{2}+2\,x
 \right) }{{{\rm e}^{2\,x}}{x}^{2}{y}^{3}+2\,{{\rm e}^{2\,x}}{x}^{4}-2
\,{{\rm e}^{x}}x{y}^{3}-4\,{{\rm e}^{x}}{x}^{3}+{y}^{3}+2\,{x}^{2}}}.
\end{equation}

\item Next we're going to apply the transformation $T_2 = \{x= \ln(x),\,y=y\}$ obtaining:
\begin{equation}
\label{1ode2ex4}
\!\!\!\!\!\!\!\!\frac{dy}{dx} \!=\! {\frac {y ( 2\,{x}^{2} \ln(x) ^{3}
+x\ln(x){y}^{3}-x \ln(x) ^{3}+x{y}^{3}-5\,x \ln(x) ^{2}
+2\,\ln(x) ) }{ ( {x}^{2} \ln(x) ^{2}{y}^{3}\!+\!2\,{x}^{2} \ln(x) ^{4}\!-\!2\,x\ln(x){y}^{3}\!-\!4\,x \ln(x) ^{3}+{y}^{3}+2\, \ln(x) ^{2} ) x}}.
\end{equation}

\item The vector field $\chi_{[2]}$ associated with the 1ODE$_{[2]}$ is defined by:
\begin{eqnarray}
\label{ex4f2}
f_2 &=& \left( {y}^{3}+2\,{z}^{2} \right)  \left( xz-1 \right) ^{2}x,
 \\ [2mm]
\label{ex4g2}
g_2 &=& y \left( 2\,{x}^{2}{z}^{3}+zx{y}^{3}+x{y}^{3}-x{z}^{3}-5\,{z}^{2}x+2\,z \right),
 \\ [2mm]
\label{ex4h2}
h_2 &=& \left( {y}^{3}+2\,{z}^{2} \right)  \left( xz-1 \right) ^{2}.
\end{eqnarray}

\item Let's apply the transformation $T_3 = \{x = y, y = z, z = x\}$ obtaining:
\begin{eqnarray}
\label{ex4f3}
f_3 &=& \left( {z}^{3}+2\,{x}^{2} \right)  \left( xy-1 \right) ^{2},
 \\ [2mm]
\label{ex4g3}
g_3 &=& \left( {z}^{3}+2\,{x}^{2} \right)  \left( xy-1 \right) ^{2}y,
 \\ [2mm]
\label{ex4h3}
h_3 &=& z \left( 2\,{x}^{3}{y}^{2}+xy{z}^{3}-{x}^{3}y+y{z}^{3}-5\,{x}^{2}y+2\,x \right).
\end{eqnarray}

\item Constructing a polynomial $P$ with arbitrary coefficients in $(x, y, z)$ and substituting in the equation (\ref{eqSsys}), we obtain (for degree 5):
\begin{equation}
\label{p23}
P = x\,z \left( -{z}^{3}+{x}^{2} \right)  \,\,\,\, \rightarrow \,\,\,\, S_{1\,[3]} = \frac{x\,z \left( -{z}^{3}+{x}^{2} \right)}{\left( {z}^{3}+2\,{x}^{2} \right)  \left( xy-1 \right) ^{2}}.
\end{equation}

\item Using the $S$-function method we can find a first integral for $\chi_{[3]}$:
\begin{equation}
\label{ii3}
I_{[3]} = -\ln  \left( {z}^{3}-{x}^{2} \right) +2\,\ln  \left( z \right) -
 \frac{1}{xy-1}.
\end{equation}

\item Applying ${T_3}^{-1}$ to $I_{[3]}$ we obtain:
\begin{equation}
\label{ii2}
I_{[2]} = -\ln  \left( {y}^{3}-{z}^{2} \right) +2\,\ln  \left( y \right) -
 \frac{1}{xz-1}.
\end{equation}

\item From $I_{[2]}$ we can write the general solution of 1ODE$_{[2]}$ ($z\,\rightarrow\,\ln(x)$):
\begin{equation}
\label{sol1ode2ex4}
-\ln  \left( {y}^{3}- \left( \ln  \left( x \right)  \right) ^{2}
 \right) +2\,\ln  \left( y \right) - \frac{1}{x\ln  \left( x \right) -1}=C_2.
\end{equation}

\item Applying ${T_2}^{-1}=\{x= {\rm e}^x,\,y=y\}$ to the general solution of 1ODE$_{[2]}$ we arrive at the general solution of 1ODE$_{[1]}$:
\begin{equation}
\label{sol1ode1ex4}
-\ln  \left( {y}^{3}-{x}^{2} \right) +2\,\ln  \left( y \right) -
 \frac{1}{x{{\rm e}^{x}}-1} = C_1.
\end{equation}

\item Finally, applying ${T_1}^{-1}=\{x= x\,y,\,y=y\}$ to the general solution of 1ODE$_{[1]}$ we have the general solution of the original 1ODE:
\begin{equation}
\label{sol1odeoex4}
-\ln  \left( -{x}^{2}{y}^{2}+{y}^{3} \right) +2\,\ln  \left( y
 \right) - \frac{1}{xy{{\rm e}^{xy}}-1} = C.
\end{equation}

\end{itemize}

\begin{obs} Some comments:
\begin{itemize}
\item The fact that we do a second transformation ($T_2$) on the transformed 1ODE (1ODE$_{[1]}$) is just a matter of standardization of the computational procedure, we could equally deal directly with the case where $\theta = {\rm e}^{x}$.

\item Regarding the $T_3$ transformation we can say something similar: It is just a matter of always working with the $S_1$ function (instead of $S_2$ and $S_3$) with the denominator being a known polynomial.

\item Although the application of a sequence of transformations followed by the sequence of inverse transformations (applied in reverse order) is a bit tiring to follow, in practice these steps have no `algorithmic weight', i.e., they consume virtually no memory or CPU time. The important point is that these transformations lead to a vector field whose associated $S$-function (to be determined) depends on obtaining only one polynomial (instead of 2). This greatly simplifies the algorithm.
\end{itemize}
\end{obs}

\subsection{Two possible algorithms}
\label{tpalg}

In this section we will present two possible algorithms to solve an 1ODE that presents an elementary function. The first is more generic, yet less efficient. The second can be more restricted, but much more efficient than the first.

\begin{algor}[$Lsolv$] This algorithm is based on the method described in the section \ref{mffi}.
\label{lsolv}

\noindent
{\bf Steps:}
\begin{enumerate}

\item Construct the vector field $\chi$, i.e., determine $[f,g,h]$ (by using definition \ref{lsset}).

\item Choose $DegMNP$ (a list of three positive integers - [$d_M$,$d_N$,$d_P$]) for the degrees of the candidates $M_c$, $N_c$ and $P_c$ (candidates for the polynomials $M_0$, $N_0$ and $P_3$).

\item Construct three polynomials $M_c$, $N_c$ and $P_c$ of degrees $d_M$, $d_N$ and $d_P$, respectively, with undetermined coefficients. 

\item Substitute them in the equation $E_1\!: M_c\,f - N_c\,h+(z\,f-g)\,P_c=0$.

\item  Collect the equation $E_1$ in the variables $(x,y,z)$ obtaining a set of (linear) equations $S_{E_1}$ for the coefficients of the polynomials $M_c$, $N_c$ and $P_c$.

\item Solve $S_{E_1}$ to the undetermined coefficients. 

\item Substitute the solution of $S_{E_1}\!\!$ in the equation $E_2\!:\! D_x(S_c) - {S_c}^2 - \partial_z( \Phi_c)\,S_c + \partial_y( \Phi_c)\!=\!0$, where $D_x \equiv \partial_x+z\,\partial_y+\Phi_c\,\partial_z$, $\Phi_c = M_c/N_c$ and $S_c = P_c/N_c$.

\item Collect the numerator of $E_2$ in the variables $(x,y,z)$ obtaining a set of equations $S_{E_2}$ for the remaining undetermined coefficients of the polynomials $M_c$, $N_c$ and $P_c$. 

\item Solve $S_{E_2}$ for the remaining undetermined coefficients. If no solution is found then FAIL.

\item Substitute the solution of $S_{E_2}$ in $M_c$ and $N_c$ to obtain $M_0$, $N_0$ and in $P_c/N_c$ to obtain $S_1$.

\item Use the $S$-function method to find the Liouvillian first integral $I$ of the 2ODE $z' = \Phi(x,y,z)$.

\item Make the substitution $z \rightarrow \theta$ on the first integral $I$ and obtain the desired general solution of the 1ODE.

\end{enumerate}

\end{algor}

\medskip

\begin{obs}
Some comments:

\begin{itemize}

\item Since there is no bound for the degrees $[d_M,d_N,d_P]$ the procedure may never get to a conclusion, i.e., more formally, $Lsolv$ is a semi algorithm.

\item Some of the steps described above are obviously much more complicated than others. Some involve very complex algorithms in themselves, for example, the algorithm that solves nonlinear polynomial systems. Besides that, since we are adapting the method of the S function to solve the final part, the penultimate item is itself a very complicated algorithm.

\end{itemize}

\end{obs}

\bigskip

\begin{algor}[$FastLs$] 
\label{fastls}
This algorithm is an improvement of the algorithm \ref{lsolv} and is based on the theorems \ref{theoAimp} and \ref{theoSsis}, on the corollary \ref{coroAimp} and on their consequences (see example 4). It is restricted to 1ODEs $\in L_S$.
\
\vspace{3mm}

{\bf Steps:}
\begin{enumerate}

\item Use a variable transformation $T_1=\{x=r(x,y),y=y\}$ where $r(x,y)$ is the argument of $\theta$ to obtain an 1ODE$_{[1]}$ such that $\theta_{[1]} = \exp(x)$ or $\theta_{[1]} = \ln(x)$.

\item If $\theta_{[1]} = \exp(x)$ then make a transformation $T_2=\{x=\ln(x),y=y\}$ to obtain an 1ODE$_{[2]}$ such that $\theta_{[2]} = \ln(x)$.

\item Construct the vector field $\chi_2=f_2\,\partial_{x}+g_2\,\partial_{y}+h_2\,\partial_{z}$ associated with the 1ODE$_{[2]}$, where 
\begin{equation}
f_2=x\,N(x,y,z),\,\,g_2=x\,M(x,y,z), \,\,h_2=N(x,y,z),
\end{equation}
if $x$ is not a factor of $N$ or 
\begin{equation}
f_2=N(x,y,z),\,\,g_2=M(x,y,z), \,\,h_2=\frac{N(x,y,z)}{x},
\end{equation}
if $x|N$.

\item Apply the transformation $T_3 = \{x = y, y = z, z = x\}$ to the vector field $\chi_{[2]}$, obtaining $\chi_3 = f_3\,\partial_{x}+g_3\,\partial_{y}+h_3\,\partial_{z}$, where
\begin{eqnarray}
f_3 &=& T_3(h_2), \nonumber 
\\ [2mm]
g_3 &=& T_3(f_2), \nonumber 
\\ [2mm]
h_3 &=& T_3(g_2). \nonumber
\end{eqnarray}

\item Choose $d_P$ for the degree of the candidate $P$ (candidate for the polynomial $P_2$ -- numerator of the $S$-function).

\item Construct a polynomial $P$ in $(x, y, z)$ with arbitrary coefficients and substitute it in the equation $E_1\!: $ (\ref{eqSsys}). 

\item  Collect the equation $E_1$ in the variables $(x,y,z)$ obtaining a set of quadratic equations $S_{E_1}$ for the coefficients of the polynomial $P$.

\item Solve $S_{E_1}$ to the undetermined coefficients. If no solution is found then FAIL.

\item Substitute the solution of $S_{E_1}\!\!$ in $P/f$ to obtain $S_{1\,[3]}$.

\item Use the $S$-function method to determine the first integral $I_{[3]}$ of the vector field $\chi_{[3]}$.

\item Apply the inverse trasformation ${T_3}^{-1}$ and obtain the first integral $I_{[2]}$ of the vector field $\chi_{[2]}$. From $I_{[2]}$ we can have (doing $z \rightarrow \ln(x)$) the general solution of the 1ODE$_{[2]}$.

\item Aplly ${T_2}^{-1}$ to the general solution of the 1ODE$_{[2]}$ and obtain the general solution of the 1ODE$_{[1]}$.

\item Aplly ${T_1}^{-1}$ to the general solution of the 1ODE$_{[1]}$ and obtain the general solution of the original 1ODE. 
\end{enumerate}

\end{algor}

\medskip

\begin{obs}
The algorithm $FastLs$ uses the equation (\ref{eqSsys}) to determine the polynomial $P$ (the numerator of the $S$-function). Since (\ref{eqSsys}) has a term of the form $\,S^2$, we have that it can be written as (see the example 4):
\begin{equation}
\label{eqSsys2}
\chi \left(\frac{P}{f}\right)=\left(\frac{P}{f}\right)^2 \frac{fg_z-gf_z}{f}+\left(\frac{P}{f}\right) \frac{gf_y-fg_y+fh_z-hf_z}{f} - \frac{fh_y-hf_y}{f}.
\end{equation}
We can isolate the equation (\ref{eqSsys2}) for $P^2$, obtaining $\,P^2={\rm polynomial}$. Since the higher degree monomials of $P$ squared cannot cancel, the maximum degree of the polynomial $P$ is bounded. In this case, therefore, the part of $FastLs$ that determines the $S$-function is a full algorithm. 
\end{obs}

\begin{obs}
The fact that $\,h | P_2,\,$ $f=x\,h\,$ and $\,g | (x\,P_1+P_3)$ (see theorem \ref{theoAimp}, case ($ii$)) does not imply that $h =c\, P_2$ and that $g = -c\,(P_1+z\,P_3)$ ($c$ constant). If the polynomials $P_2$ and $(x\,P_1+P_3)$ have a polynomial factor in common, the $FastLs$ algorithm may lose this case. 
\end{obs}

\bigskip

\section{Performance}
\label{performance}
\hspace\parindent
The analysis presented in this section is divided in two parts\footnote{In this paper all the computational data (time of running etc) was obtained on the same computer with the following configuration: Intel(R) Core(TM) i5-8265U @ 1.8 GHz.}:
\begin{itemize}
\item In the first subsection, we solve (with our package -- see appendix) a set of LLS equations that we could not solve using the traditional methods (and even non-traditional methods) implemented in the Maple computer algebra system (Maple CAS).
\item In the second, we present an example of a LLS equation with parameters and use the integrability analysis (provided by using the {\tt Par} parameter in the {\tt Sfunction} command -- see appendix) to find values of these parameters that allow the existence of a Liouvillian first integral.
\end{itemize}

\subsection{A set of `hard' LLS equations}
\label{hard1odes}
\hspace\parindent
In this section we present a set of LLS equations which the powerful methods implemented in the Maple ODE solver (release 17) -- {\tt dsolve} command -- cannot solve. Bellow, we show the CPU time and memory that our package ($LeapsLS$, that implements our procedure -- see appendix) spends to find the solutions (through the use of {\tt Solde} command). In order to avoid presenting more obvious steps, we will start by presenting the 1ODEs resulting from the use of the transformation (\ref{tlevsmi}) in the LLS equations.

\medskip

Consider the following 1ODEs:
\begin{enumerate}

\item  The following five 1ODEs have an elementary Liouvillian solution. We will refer to them as 1ODEs 1, 2, 3, 4, 5:
\begin{equation}
\label{harde1}
\frac{dy}{dx}={\frac {x \left( 2\, \left( \ln  \left( x \right)  \right) ^{2}{x}^{4}
{y}^{2}+2\,\ln  \left( x \right) {x}^{2}y+2\,\ln  \left( x \right) {y}
^{2}-y{x}^{2}+{y}^{2}+2 \right) }{ \left( \ln  \left( x \right) 
 \right) ^{2}{x}^{4}{y}^{2}+\ln  \left( x \right) {x}^{4}+\ln  \left( 
x \right) {x}^{2}y+1}}
\end{equation}

\begin{equation}
\label{harde2}
\frac{dy}{dx}={\frac {3 \left( {{\rm e}^{x}} \right) ^{2}\!{x}^{4}\!{y}^{4}\!+\! \left( {y}^{4}\!-\!{
y}^{3}\!{x}^{4}\!-\!{x}^{3}\!{y}^{3}\!+\!6\,{x}^{3}\!{y}^{2}+x{y}^{4}\!
 \right) {{\rm e}^{x}}+3\,{x}^{2}}{ \left( {{\rm e}^{x}} \right) ^{2}{
x}^{2}{y}^{4}+ \left( {x}^{4}{y}^{2}-x{y}^{3}+2\,x{y}^{2} \right) {
{\rm e}^{x}}-{x}^{3}+y+1}}
\end{equation}

\begin{equation}
\label{harde3}
\frac{dy}{dx}=-{\frac { \left( 2\,{x}^{3}{y}^{4}-4\,{x}^{2}{y}^{2}-{x}^{2}+\ln 
 \left( y \right) +2\,x \right) y}{-2\,{x}^{4}{y}^{2}-{x}^{2}{y}^{4}+2
\,\ln  \left( y \right) {x}^{2}{y}^{2}+2\,x{y}^{2}-1}}
\end{equation}

\begin{equation}
\label{harde4}
\frac{dy}{dx}=-{\frac {3\, \left( {{\rm e}^{y}} \right) ^{2}{x}^{4}{y}^{2}-7\,{
{\rm e}^{y}}{x}^{3}y+3\,{x}^{2}-y}{x \left(  \left( {{\rm e}^{y}}
 \right) ^{2}{x}^{4}{y}^{2}-3\,{{\rm e}^{y}}{x}^{3}y-{x}^{3}{{\rm e}^{
y}}+{x}^{2}-y-1 \right) }}
\end{equation}

\begin{equation}
\label{harde5}
\frac{dy}{dx}=-{\frac {4\,{x}^{3}{y}^{6}+8\,{x}^{4}{y}^{3}+4\,{x}^{5}-{x}^{4}+\cos
 \left( y \right) }{ \left( {y}^{6}+2\,x{y}^{3}+{x}^{2} \right) \sin
 \left( y \right) -3\,{x}^{4}{y}^{2}+3\,\cos \left( y \right) {y}^{2}}}
\end{equation}

\bigskip

\item  The following five 1ODEs have a non-elementary Liouvillian solution. We will refer to them as 1ODEs 6, 7, 8, 9, 10. 
\begin{equation}
\label{harde6}
\frac{dy}{dx}={\frac { \left( 2\,{x}^{2}-y \right)  \left( \ln  \left( x \right) {y}
^{2}-{x}^{2}y-1 \right) y}{x \left(  \left( \ln  \left( x \right) 
 \right) ^{2}{y}^{3}-\ln  \left( x \right) {x}^{2}{y}^{2}- \left( \ln 
 \left( x \right)  \right) ^{2}{y}^{2}+2\,\ln  \left( x \right) {x}^{2
}y-{x}^{4}-\ln  \left( x \right) y \right) }}
\end{equation}

\begin{equation}
\label{harde7}
\frac{dy}{dx}=-{\frac { \left( \ln  \left( x \right)  \right) ^{2}x{y}^{2}-2\,\ln 
 \left( x \right) {x}^{2}y-\ln  \left( x \right) xy+\ln  \left( x
 \right) {y}^{2}+{x}^{3}}{x \left( \ln  \left( x \right)  \right) ^{2}y}}
\end{equation}

\begin{equation}
\label{harde8}
\frac{dy}{dx}={\frac { \left( -{{\rm e}^{2\,x}}{x}^{6}{y}^{2}-2\,{{\rm e}^{2\,x}}{x}
^{5}{y}^{2}+{{\rm e}^{2\,x}}{x}^{4}{y}^{2}+{{\rm e}^{x}}{x}^{4}y+3\,{
{\rm e}^{x}}{x}^{3}y+1 \right) {{\rm e}^{-x}}}{{x}^{3} \left( {{\rm e}
^{x}}{x}^{3}y-x-1 \right) }}
\end{equation}

\begin{equation}
\label{harde9}
\frac{dy}{dx}=-{\frac { \left( \ln  \left( y \right)  \right) ^{2}x{y}^{4}- \left( 
\ln  \left( y \right)  \right) ^{2}{y}^{4}-{x}^{2}\ln  \left( y
 \right) {y}^{2}+\ln  \left( y \right) x{y}^{2}+{x}^{3}}{xy \left( 2\,
\ln  \left( y \right) +1 \right)  \left( \ln  \left( y \right) {y}^{2}
-{x}^{2}-x \right) }}
\end{equation}

\begin{equation}
\label{harde10}
\frac{dy}{dx}={\frac {{y}^{3} \left( -x{y}^{3}+{{\rm e}^{y}}+1 \right) }{{{\rm e}^{y
}}{x}^{2}{y}^{6}+3\,{x}^{2}{y}^{5}-2\,{{\rm e}^{y}}x{y}^{3}-3\,{
{\rm e}^{y}}x{y}^{2}-3\,x{y}^{2}+{{\rm e}^{y}}}}
\end{equation}

\end{enumerate}

\bigskip

Applying the {\tt Solde} command to each one of them, we have:
\begin{verbatim}
[> t0 := time(): solde[1] := Solde(_1ode[1],Deg=7); time()-t0;
\end{verbatim}
$$
solde[1] :=-{\frac {\ln  \left( {x}^{2}-y \right) \ln  \left( x \right) {x}^{2}y+
\ln  \left( {x}^{2}-y \right) +1}{\ln  \left( x \right) {x}^{2}y+1}}={
\it \_C1}
$$
$$ 0.563 $$

\medskip

\begin{verbatim}
[> t0 := time(): solde[2] := Solde(_1ode[2],Deg=7); time()-t0;
\end{verbatim}
$$
solde[2] := -{\frac {x{y}^{2}{{\rm e}^{x}}+1}{{{\rm e}^{x}}\ln  \left( {x}^{3}-y
 \right) x{y}^{2}+\ln  \left( {x}^{3}-y \right) +y}}={\it \_C1}
$$
$$ 1.078 $$

\medskip

\begin{verbatim}
[> t0 := time(): solde[3] := Solde(_1ode[3],Deg=5); time()-t0;
\end{verbatim}
$$
solde[3] := \left( -{x}^{2}+\ln  \left( y \right)  \right) {{\rm e}^{{\frac {x}{x
{y}^{2}-1}}}}={\it \_C1}
$$
$$ 0.209 $$

\medskip

\begin{verbatim}
[> t0 := time(): solde[4] := Solde(_1ode[4],Deg=5); time()-t0;
\end{verbatim}
$$
solde[4] := -{\frac {{{\rm e}^{y}}xy\ln  \left( {x}^{3}{{\rm e}^{y}}+1 \right) -
\ln  \left( {x}^{3}{{\rm e}^{y}}+1 \right) +1}{xy{{\rm e}^{y}}-1}}={
\it \_C1}
$$
$$ 0.625 $$

\medskip

\begin{verbatim}
[> t0 := time(): solde[5] := Solde(_1ode[5],Deg=7); time()-t0;
\end{verbatim}
$$
solde[5] := \ln  \left( -2\,{{\rm e}^{iy}}{x}^{4}+{{\rm e}^{2\,iy}}+1 \right) +{
\frac {-i{y}^{4}-iyx+1}{{y}^{3}+x}}={\it \_C1}
$$
$$ 3.687 $$

\medskip

\begin{verbatim}
[> t0 := time(): solde[6] := Solde(_1ode[6],Deg=5); time()-t0;
\end{verbatim}
$$
solde[6] := - \frac{{\it Ei} \left( 1,- \left( \ln  \left( x \right) y-{x}^{2}
 \right) ^{-1} \right) y+{{\rm e}^{ \left( \ln  \left( x \right) y-{x}
^{2} \right) ^{-1}}}}{y}={\it \_C1}
$$
$$ 0.203 $$

\medskip

\begin{verbatim}
[> t0 := time(): solde[7] := Solde(_1ode[7],Deg=7); time()-t0;
\end{verbatim}
$$
solde[7] := \left( {\it Ei} \left( 1, \left( \ln  \left( x \right) y-x \right) ^{
-1} \right) {{\rm e}^{ \left( \ln  \left( x \right) y-x \right) ^{-1}}
}+x \right) {{\rm e}^{- \left( \ln  \left( x \right) y-x \right) ^{-1}
}}={\it \_C1}
$$
$$ 0.218 $$

\medskip

\begin{verbatim}
[> t0 := time(): solde[8] := Solde(_1ode[8],Deg=9); time()-t0;
\end{verbatim}
$$
solde[8] := \left( {\it Ei} \left( 1,- \left( {{\rm e}^{x}}{x}^{2}y-1 \right) ^{-
1} \right) x+{{\rm e}^{ \left( {{\rm e}^{x}}{x}^{2}y-1 \right) ^{-1}}}
 \right) {x}^{-1}={\it \_C1}
$$
$$ 0.687 $$

\medskip

\begin{verbatim}
[> t0 := time(): solde[9] := Solde(_1ode[9],Deg=6); time()-t0;
\end{verbatim}
$$
solde[9] := x\,{{\rm e}^{\displaystyle{\frac {x}{\ln  \left( y \right) {y}^{2}-x}}}}+{\it Ei}
 \left( 1,-{\frac {x}{\ln  \left( y \right) {y}^{2}-x}} \right) ={\it 
\_C1}
$$
$$ 0.437 $$

\medskip

\begin{verbatim}
[> t0 := time(): solde[10] := Solde(_1ode[10],Deg=6); time()-t0;
\end{verbatim}
$$
\left( {{\rm e}^{- \left( x{y}^{3}-1 \right) ^{-1}}}{\it Ei} \left( 1
,- \left( x{y}^{3}-1 \right) ^{-1} \right) +{{\rm e}^{y}} \right) {
{\rm e}^{ \left( x{y}^{3}-1 \right) ^{-1}}}={\it \_C1}
$$
$$ 0.297 $$

\medskip

\begin{tabular}
{|c|c|c|}
\hline
 1ODE & Time (sec.) & Memory (MB) \\
\hline
 1 (\ref{harde1}) & 0,672 & 21,07\\ 
\hline
 2 (\ref{harde2}) & 1,078 & 35,60\\ 
\hline
 3 (\ref{harde3}) & 0,209 & 4,18\\ 
\hline
 4 (\ref{harde4}) & 0,625 & 33,18\\ 
\hline
 5 (\ref{harde5}) & 3,687 & 28,47\\ 
\hline
 6 (\ref{harde6}) & 0,329 & 17,41\\ 
\hline
 7 (\ref{harde7}) & 0,218 & 16,18\\ 
\hline 
 8 (\ref{harde8}) & 0,687 & 21,26\\ 
\hline
 9 (\ref{harde9}) & 0,437 & 29,86\\ 
\hline
 10 (\ref{harde10}) & 0,297 & 18,42\\
\hline
\end{tabular}

\medskip

\begin{obs}
\label{hardescomm}
Some comments:
\begin{itemize}
    \item When we start a Maple session, the simple loading of the basic packages that are activated with the opening of the program's standard platform already uses 4.18 MB. So, when that number appears (in the table above), the meaning is basically `almost no memory expenditure'.

    \item For the measurement of consumed CPU time, a similar problem occurs: the time measurement provided is `variable from session to session'. However, the number representing the time spent is always reasonably close to the average of the values obtained in different rounds. In the data placed in the table above, we avoid averaging over many sessions (which means that the last two digits are not significant) because we were focused on a more qualitative result.
    
    \item Another piece of information (quite relevant for those who don't frequently use the Maple command {\tt dsolve}) is that the memory consumed by our package in each one of the ten 1ODEs is relatively low.
    
    \item Another interesting detail is that most of the time consumed (even in cases where the time consumed is very short) is generally used in determining the $S$-function, that is, in addition to validating our method, the previous result is also a reaffirmation of the $S$-function method.

    \item Finally, looking at the table below, we can see why these 1ODEs cause difficulties for the vast majority of methods. We only need to pay attention to the Lie symmetries or to the integrating factors of the 1ODEs' set. We see that only 1ODEs 3, 5 and 10 have integrating factors in which the $\theta$ function is absent (and, even so, they are not simple).
\end{itemize}
\end{obs}

\medskip

\begin{tabular}
{|c|c|c|}
\hline
 1ODE & Symmetry ($\nu$) & Integrating Factor \\
\hline
 1 (\ref{harde1}) & $-{\frac { ( {x}^{2}y\ln  \left( x \right) +1 ) ^{2}}{{
{\rm e}^{ \left( {x}^{2}y\ln  \left( x \right) +1 \right) ^{-1}}}
 \left(  \left( \ln  \left( x \right)  \right) ^{2}{x}^{4}{y}^{2}+\ln 
 \left( x \right) {x}^{4}+{x}^{2}y\ln  \left( x \right) +1 \right) }}
$ & ${\frac {{{\rm e}^{ ( {x}^{2}y\ln  \left( x \right) +1 ) ^{-1}}}}{ \left( {x}^{2}y\ln  \left( x \right) +1 \right) ^{2}}}$\\ 
\hline
 2 (\ref{harde2}) & $- \frac{ ( x{y}^{2}{{\rm e}^{x}}+1 ) ^{2} \, {{\rm e}^{{\frac 
{-y}{x{y}^{2}{{\rm e}^{x}}+1}}}}  }{  \left( {{\rm e}^
{x}} \right) ^{2}{x}^{2}{y}^{4}+{{\rm e}^{x}}{x}^{4}{y}^{2}-{{\rm e}^{
x}}x{y}^{3}+2\,x{y}^{2}{{\rm e}^{x}}-{x}^{3}+y+1 }
$ & $\frac{{\rm e}^{{\frac {y}{x{y}^{2}{{\rm e}^{x}}+1}}}}{ \left( x{y}^{2}{{\rm e}^{x}}+1 \right) ^{2}}$ \\ 
\hline
 3 (\ref{harde3}) & $ \frac{ ( x{y}^{2}-1 ) ^{2}y \, {{\rm e}^{{\frac {-x}{x{y}^{2}-
1}}}} }{ -2\,{x}^{4}{y}^{2}-{x}^{2}{y}^{4}+2\,\ln 
 \left( y \right) {x}^{2}{y}^{2}+2\,x{y}^{2}-1 }
$ & $- \frac{{\rm e}^{{\frac {x}{x{y}^{2}-1}}}}{y\, \left( x{y}^{2}-1 \right) ^{2}}$ \\ 
\hline
 4 (\ref{harde4}) & ${\frac { \left( xy{{\rm e}^{y}}-1 \right) ^{2}\,{{\rm e}^{ \left(1- xy{
{\rm e}^{y}} \right) ^{-1}}}}{{{\rm e}^{y}}x \left(  \left( {{\rm e}^{
y}} \right) ^{2}{x}^{4}{y}^{2}-3\,{{\rm e}^{y}}{x}^{3}y-{x}^{3}{
{\rm e}^{y}}+{x}^{2}-y-1 \right) }}
$ & ${\frac {-{{\rm e}^{ \left( xy{{\rm e}^{y}}-1 \right) ^{-1}}}{{\rm e}^{y}}}{ \left( xy{{\rm e}^{y}}-1 \right) ^{2}}}$ \\ 
\hline
 5 (\ref{harde5}) & ${\frac { ( {y}^{3}+x )^{2}\,{{\rm e}^{ ( -{y}^{3}-x
) ^{-1}}}}{ \left( \sin \left( y \right) {y}^{6}-3\,{x}^{4}{y}^{2
}+2\,\sin \left( y \right) x{y}^{3}+3\,\cos \left( y \right) {y}^{2}+
\sin \left( y \right) {x}^{2} \right) }}
$ & ${\frac {-{{\rm e}^{ ( {y}^{3}+x ) ^{-1}}}}{ \left( {y}^{3}+x \right) ^{2}}}$ \\ 
\hline
 6 (\ref{harde6}) & ${\frac { ( \ln  \left( x \right) y-{x}^{2} ) ^{2}{y}^{2}\,{{\rm e}^{({x}^{2} -  \ln  \left( x \right) y) ^{-1}}}}{
  \left( \ln  \left( x \right)  \right) ^{2}{y}^{3}-\ln  \left( x \right) {x}^{2}{y}^{2}- \left( \ln  \left( x \right) 
 \right) ^{2}{y}^{2}+2\,\ln  \left( x \right) {x}^{2}y-{x}^{4}-\ln 
 \left( x \right) y  }}
$ & ${\frac {-{{\rm e}^{ ( \ln  \left( x \right) y-{x}^{2}) ^{-1}}}}{ \left( \ln  \left( x \right) y-{x}^{2} \right) ^{2}x{y}^{2}}}$ \\ 
\hline
 7 (\ref{harde7}) & $\frac{ ( \ln  \left( x \right) y-x ) ^{2} }{ {{\rm e}^{-
 \left( \ln  \left( x \right) y-x \right) ^{-1}}} \, ( \ln  \left( x \right)  ) ^{2}\,y}
$ & $\frac{-{\rm e}^{- \left( \ln  \left( x \right) y-x \right) ^{-1}}}{( 
\ln  \left( x \right) y-x ) ^{2}\,x}$ \\ 
\hline 
 8 (\ref{harde8}) & $\frac{ {{\rm e}^{-{\frac {{{\rm e}^{x}}{x}^{3}y-x+1}{{{\rm e}^{x}}{x}^{2}y-1}
}}} ( {{\rm e}^{x}}{x}^{2}y-1 ) ^{2} }{ {x}\, \left( {
{\rm e}^{x}}{x}^{3}y-x-1 \right) }
$ & $\frac{-{{\rm e}^{{\frac {{{\rm e}^{x}}{x}^{3}y-x+1}{{{\rm e}^{x}}{x}^{2}y-1}
}}} }{ {x}^{2} \left( {{\rm e}^{x}}{x}^{2}y-1 \right) ^{2}}$ \\ 
\hline
 9 (\ref{harde9}) & $ \frac{ {{\rm e}^{-{\frac {x}{\ln  \left( y \right) {y}^{2}-x}}}} ( \ln 
 \left( y \right) {y}^{2}-x ) ^{2} }{ {y}\, \left( 2\,\ln 
 \left( y \right) +1 \right)\,\left( \ln  \left( y \right) {y}^{2
}-{x}^{2}-x \right) }
$ & $ \frac{ -{{\rm e}^{{\frac {x}{\ln  \left( y \right) {y}^{2}-x}}}} }{ \left( \ln  \left( y \right) {y}^{2}-x \right) ^{2}\,{x} }$ \\ 
\hline
 10 (\ref{harde10}) & $\frac{ {{\rm e}^{- ( x{y}^{3}-1 ) ^{-1}}} ( x{y}^{3}-1 ) ^{2} }{ {{\rm e}^{y}}{x}^{2}{y}^{6}+3\,{x}^{2}{y}^{5}-2\,
{{\rm e}^{y}}x{y}^{3}-3\,{{\rm e}^{y}}x{y}^{2}-3\,x{y}^{2}+{{\rm e}^{y}}  }
$ & ${\frac {-{{\rm e}^{ ( x{y}^{3}-1 ) ^{-1}}}}{ \left( x{y}^{3}-1 \right) ^{2}}}$ \\
\hline
\end{tabular}

\begin{obs}
For all 1ODEs in Kamke's book \cite{Kam} that obey the conditions of application of the algorithm our method does very well in all cases. As they are much simpler cases than those shown in the table, we decided not to include them in order not to increase the presentation unnecessarily.
\end{obs}

\subsection{Integrability analisys}
\label{ub}
\hspace\parindent
In this section, we will present a very useful way to apply the $ LeapsLS $ package: the {\tt Par} parameter (see appendix). When we are studying an ODE presenting parameters, we can use the parameter {\tt Par} in the entries of the {\tt Sfunction} command. Thus, if the command was unable to find a $S$-function for generic parameter values, the entry {\tt Par = \{$\alpha $,$\beta$, ...\}} tells the command to try to find a $S$-function for some relationship between them, i.e., the method brings us the possibility of an analysis of the integrability region of the 1ODE's parameters\footnote{This happens because of the algebraic nature of the method used to determine the coefficients of the polynomial that represents the numerator of the $S$ function: we can add the parameters present in the ODE to the undetermined coefficients of the $P_2$ polynomial ($S$ numerator) and the solution, if any, emerges directly from the process.}. Let's analyze the following case:
\begin{equation}
\label{levsmi}
\ddot{x} + (d\,{x}^{2}{{\rm e}^{x}}+b\,x{{\rm e}^{x}}+c\,{y}+a)\,\dot{x} + \alpha\,{x}^{2}{{\rm e}^{x}}+\beta\,x = 0,
\end{equation}
Applying the transformation (\ref{tlevsmi}) we obtain the 1ODE 
\begin{equation}
\label{levsminos}
\frac{dy}{dx}= {\frac {(-d{x}^{2}{{\rm e}^{x}}-b\,x{{\rm e}^{x}}-c\,{y}-a)\,y-\alpha\,{x}^{2}{{\rm e}^{x}}-\beta\,x}{y}}
\end{equation}  
Applying the {\tt Sfunction} command and using the parameter {\tt Par} (see appendix)
\begin{verbatim}
[> _1ode := diff(y(x),x)= \cdots:
[> ss1 := [Sfunction(_1ode,Par={b,d,e},DegMNP=[6,4,5])]:
[> for i from 3 by 3 to nops(ss1) do 
   for j to nops(ss1[i]) do 
   print([i,j],___________________,traperror(factor(subs(ss1[i][j],
   ss1[i-2]))),subs(ss1[i][j],[`b=`,b,`c=`,c,`B=`,B,`A=`,A,`d=`,d]));
   trys[i,j]:=traperror(factor(subs(ss1[i][j],ss1[i-1][1]/ss1[i-1][2])));
   print(trys[i,j]); 
   end do; 
   end do;
\end{verbatim}
we obtain\footnote{We omit the full output because it is too long.}
$$\vdots $$
$$ [3, 4], \_\_\_\_\_\_\_\_ \,, -{\frac {cz}{cy+a}}\, [`b= `, b, `c= `, c, `B= `, 0, `A= `, 0, `d= `, d] $$
$$ {\frac { \left( c d{x}^{2}z+bcxz+{c}^{2}y+ac+cy+cz+a \right) z}{cy+a}} $$
$$\vdots $$
Using the determined region on 1ODE (\ref{levsminos}) we finally obtain
\begin{verbatim}
[> new1ode := subs(ss1[3][4],_1ode);
\end{verbatim}
\begin{equation}
\label{newde}
\,\,\,\,\,\,\,\,\,\,\,\,\,\,\,\,\,\,\,\,\,{\frac {d}{dx}}y \left( x \right) ={\frac {-d{x}^{2}y \left( x
 \right) {{\rm e}^{x}}-x{{\rm e}^{x}}y \left( x \right) b- \left( y
 \left( x \right)  \right) ^{2}c-y \left( x \right) a}{y \left( x
 \right) }}
\end{equation}
that can be solved with the {\tt Solde} command
\begin{verbatim}
[> s1 := ss1[1];
[> Solde(new1ode,SF=s1):
\end{verbatim}
$$ s1 := -{\frac {cz}{cy+a}} $$

\noindent
Rearranging the solution we have
$$
y  = \left( -\frac {d\,{\rm e}^x \left( \left( c+1 \right) ^{2}{x}^{2}-2\, \left( c+1 \right) x+2 \right) }{ \left( c+1 \right) ^3}-
{\frac {b\,{\rm e}^x \left(  \left( c+1 \right) x-1 \right) }{ \left( c+1 \right) ^{2}}}-{\frac {a}{c}}+{\it \_C_1} \right) 
$$

\noindent
Replacing $ y $ with $ \dot{x} $ and solving for ${\it \_C_1}$, we obtain the first integral of the LLS equation (for the parameter regions we find):
$$ I = \dot{x}+\left( \frac {d\,{\rm e}^x \left( \left( c+1 \right) ^{2}{x}^{2}-2\, \left( c+1 \right) x+2 \right) }{ \left( c+1 \right) ^3}+
{\frac {b\,{\rm e}^x \left(  \left( c+1 \right) y-1 \right) }{ \left( c+1 \right) ^{2}}}+{\frac {a}{c}}\right) 
$$

\section{Conclusion}
\label{conclu}
\hspace\parindent

In this work we developed a new method to compute Liouvillian first integrals of LLS equations presenting an elementary function. The method is based on a connection between a LLS equation with an elementary function $\theta$ and a polynomial vector field in three variables. The main point is that, if this LLS equation has a Liouvillian first integral, we can use the $S$-function method (see \cite{Noscpc2019}) to find a first integral for the vector field and (using the simple substitution $z \rightarrow \theta$) find the first integral of the LLS equation. The method proved to be effective even in cases where the search for first integrals by procedures such as the Lie symmetries method and Darboux approaches are particularly problematic.

The method we developed first ($ Lsolv $ algorithm) was (in a first assessment) reasonably efficient compared to the more canonical alternative that would be to try a DPS approach or the Lie symmetry method. However, the appearance of an unexpected result, improved the efficiency of our algorithm far surpassing our expectations. These result (expressed in the theorems \ref{theoAimp} and \ref{theoSsis}) allowed the creation of the $ FastL_S $ algorithm that, for a very comprehensive range of 1EDOs (the class $ L_S $ - see the definition \ref{lsset}), performed very well.

The algorithms $ Lsolv $ and $ FastL_S $ were implemented in a computational package ($LeapsLS$) on the Maple symbolic computing platform. The package, in addition to finding the first integrals using the methods mentioned above, presents commands (see appendix) that materialize all stages of the process. The commands also make use of several parameters that help both in the search for solutions as well as in the research in physics and mathematics.

Finally, due to the nature of the algorithm search process, we can perform an analysis of integrability regions LLS equations that present parameters (see the example in section \ref{ub}).


\newpage

\appendix{The $LeapsLS$ Package}

\bigskip
Here we present a Maple (release 17) implementation of the algorithms \ref{lsolv} and \ref{fastls} (with some modifications in order to improve flexibility and practicality):
\begin{itemize}
\item Title: $LeapsLS$ -- Liouvillian (elementary, algebraic, polynomial) solutions of the Liénard System.
\item All the computational data (time of running etc) was obtained on the same computer with the following configuration: 
Intel(R) Core(TM) i5-8265U @ 1.8 GHz. 
\item The operating system under which the program has been tested was Windows 10. 
\item Our implementation can find solutions in many cases where the LLS equation under study can not be solved by other powerful solvers nor by other powerful methods. 
\item The program can perform (in a very natural way) an analysis of the integrability region of the LLS equation's parameters. 
\item The package presents some useful research commands.
\end{itemize}

\section{Summary of the commands:}

\begin{itemize}
\item {\tt Lls1ode} constructs the 1ODE (\ref{1levsmi}) from the LLS equation (\ref{levsmi}).

\item {\tt Dx} constructs the $D_x$ operator (total derivative over the solutions) associated with the 1ODE.

\item {\tt Xi} constructs the $\chi$ operator associated with the 1ODE.

\item {\tt Tr1ode} applies a transformation to the 1ODE and returns the transformed 1ODE and the inverse transformation.

\item {\tt Ode2} tries to determine the 2ODE associated with the vector field $\chi$.

\item {\tt Sfunction} tries to determine a $S$-function associated with the vector field $\chi$.

\item {\tt Ode1a} determines the rational 1ODE associated with the 2ODE (from the knowledge of a $S$-function).

\item {\tt Hfunction} tries to find the general solution for an associated rational 1ODE.

\item {\tt PDEassol} constructs and tries to solve the 1PDEs that relate the first integral $I$ with the $H$-functions.

\item {\tt Solde} tries to determine a Liouvillian first integral of the LLS equation.

\end{itemize}

\section {Package commands}
\label{packcom}

Here we present a description of the main commands of the $LeapsLS$ package.

\noindent
\subsection{Command: {\tt Dx}}
\label{dxcom}

\noindent {\it Feature:} This command returns the operator $D_x$ (total derivative over the solutions).

\bigskip

\noindent
{\it Calling sequence:}

\begin{verbatim}
[> Dx(ode);
\end{verbatim}

\noindent
{\it Parameters:}
\medskip

\hspace
\parindent
{\tt ode} - The 1ODE.
\medskip

\bigskip

\noindent
{\it Synopsis:}
\smallskip
\smallskip

\noindent
The command {\tt Dx} returns the {\it total derivative operator} over the solutions of the 1ODE, i.e., the operator $d_x=\partial_x+y'\,\partial_y$ restricted to bound $I(x,y)=c$ (the solution-curves). Since, over the solutions $y'=\phi(x,y)$ this implies that 
\begin{equation}
D_x= (\partial_x+y'\,\partial_y)|_{y'=\phi} = \partial_x+\phi\,\partial_y, \,\,\,\,\left(y'=\frac{dy}{dx}\right).
\end{equation}

\bigskip


%%%%%%%
\noindent
\subsection{Command: {\tt Xi}}
\label{xicom}

\noindent {\it Feature:} This command returns the vector field $\chi$.

\bigskip

\noindent
{\it Calling sequence:}

\begin{verbatim}
[> Xi(ode);
\end{verbatim}

\noindent
{\it Parameters:}
\medskip

\hspace
\parindent
{\tt ode} - The 1ODE.
\bigskip

\noindent
{\it Extra parameters:}
\medskip

\hspace
\parindent
{\tt FGH} - If this parameter is present in the command call then the command returns a list: $[f,g,h]$.
\bigskip

\noindent
{\it Synopsis:}
\smallskip
\smallskip

\noindent
The command {\tt Xi} returns the {\it vector field $\chi$ associated with the 1ODE}, i.e., the command returns an operator $u\rightarrow f\,\partial_x(u)+g\,\partial_y(u)+h\,\partial_z(u)$, where $f$, $g$ and $h$ are polynomials given by (\ref{deefe},\ref{dege},\ref{deaga}). If the extra parameter {\tt FGH} is present in the entries then the command returns a list with the polynomials $f$, $g$ and $h$: $[f,g,h]$.

\bigskip


%%%%%%%
\noindent
\subsection{Command: {\tt Tr1ode}}
\label{tr1odecom}

\noindent {\it Feature:} This command returns the transformed 1ODE and the inverse transformation.

\bigskip

\noindent
{\it Calling sequence:}

\begin{verbatim}
[> Tr1ode(ode);
\end{verbatim}

\noindent
{\it Parameters:}
\medskip

\hspace
\parindent
{\tt ode} - The 1ODE.
\bigskip

\noindent
{\it Extra parameters:}
\medskip

\hspace
\parindent
{\tt TR=[x=F(x,y),y=G(x,y)]} - Where {\tt [x=F(x,y),y=G(x,y)]} is an invertible transformation. The default is the identity transformation.
\bigskip

\noindent
{\it Synopsis:}
\smallskip
\smallskip

\noindent
The command {\tt Tr1ode} returns a list with two entries: the {\it transformed 1ODE} and the inverse transformation of {\tt TR}. If the inverse transformation cannot be obtained, the command returns only the 1ODE transformed. If no transformation is provided, the default is the identity transformation.
\bigskip

\noindent
\subsection{Command: {\tt Ode2}}
\label{ode2com}

\noindent {\it Feature:} This command returns the 2ODE(s) associated with the 1ODE.

\bigskip

\noindent
{\it Calling sequence:}

\begin{verbatim}
[> Ode2(ode);
\end{verbatim}

\noindent
{\it Parameters:}
\medskip

\hspace
\parindent
{\tt ode} - The 1ODE.
\bigskip

\noindent
{\it Extra parameters:}
\medskip

\hspace
\parindent
{\tt PS} - If this parameter is present in the command call, then the command returns the polynomials $P_1$, $P_2$ and $P_3$: $[P_1,P_2,P_3]$. 
\bigskip


\noindent
{\it Synopsis:}
\smallskip
\smallskip

\noindent
The command {\tt Ode2} returns the function $\Phi(x,y,z) = M_0(x,y,z)/N_0(x,y,z)$ of the 2ODE associated with the 1ODE. This 2ODE is one of six possible 2ODEs, namely, the one that considers $x$ as the independent variable, $y$ as the dependent variable and $z$ as the first derivative. If the extra parameter {\tt PS} is present in the entries then the command returns a list with the polynomials $P_1$, $P_2$ and $P_3$: $[P_1,P_2,P_3]$.
\bigskip


%%%%%%%
\noindent
\subsection{Command: {\tt Sfunction}}
\label{sfuncom}

\noindent {\it Feature:} This command tries to find the $S$-function $S_1$.
\bigskip

\noindent
{\it Calling sequence:}

\begin{verbatim}
[> Sfunction(ode);
\end{verbatim}

\noindent
{\it Parameters:}
\medskip

\hspace
\parindent
{\tt ode} - The 1ODE.
\bigskip

\noindent
{\it Extra parameters:}
\medskip

\hspace
\parindent
{\tt Deg = n} - Where {\tt n} is a positive integer denoting the degree of the polynomial $P$ (the numerator or denominator of $S_1$). The default is 3.
\medskip

\hspace
\parindent
{\tt DegMNP = [degM,degN,degP]} - Where {\tt [degM,degN,degP]} is a list with the degrees of the candidates for $M_0$, $N_0$ and $P_2$, respectively.
\medskip

\hspace
\parindent
{\tt Par = $\{\alpha,\,\beta, \cdots \,\}$} - Where $\{\alpha,\,\beta, \cdots \,\}$ is a set with parameters present in the 1ODE.
\medskip

\hspace
\parindent
{\tt Den = deno} - Where {\tt deno} is the denominator (presumed) of the $S$-function.
\bigskip

\noindent
{\it Synopsis:}
\smallskip
\smallskip

The command {\tt Sfunction} tries to find a $S$-function associated with the rational 2ODE ($\Phi = M_0/N_0$) through a Liouvillian first integral $I$. The command computs (if possible) the rational function $S_1$ or a polynomial $P$ that is the numerator (or denominator) of $S_1$ (this case only if we are using the $FastLs$ algorithm). The $S$-function is the basis for finding the Liouvillian first integral of the vector field $\chi$ associated with the 1ODE (see section \ref{1odewfpvf}). In the case where the numerator or the denominator of $S_1$ is (presumed) one of the polynomials $\{f,g,h\}$, we can use the parameter {\tt Deg} to inform the program of the degree we will use for the candidate $P_c$. Associated with this case we have also the parameter {\tt Den} that can acelerate the process by indicating the numerator (or denominator -- this case only if we are using the $FastLs$ algorithm) of $S_1$. In the case where none of the polynomials $\{f,g,h\}$ is the numerator (or denominator) of $S_1$, we use the parameter {\tt DegMNP} which is, actually, a list with three integers {\tt [degM,degN,degP]} that tell the command the degree of candidates $M_c$, $N_c$ and $P_c$, respectively. Furthermore, when we are dealing with an 1ODE that has physical constants (undefined), we can use the {\tt Par} parameter to do an (Liouvillian) integrability analysis in case the command does not find a $S$-function for arbitrary values of these constants.

\bigskip

\noindent
\subsection{Command: {\tt Ode1a}}
\label{exodcom}

\noindent {\it Feature:} This command determines the associated 1ODE.
\bigskip

\noindent
{\it Calling sequence:}

\begin{verbatim}
[> Ode1a(ode);
\end{verbatim}

\noindent
{\it Parameters:} The same as above.
\medskip

\noindent
{\it Extra parameters:} The same as above and:
\medskip

\hspace
\parindent
{\tt SF = S1} - Where {\tt S1} is a $S$-function $S_1$.
\bigskip

\noindent
{\it Synopsis:}
\smallskip
\smallskip

The command {\tt Ode1a} uses the $S$-function to construct the associated 1ODE. The first extra parameters are the same of the {\tt Sfunction} command and the last extra parameter ({\tt SF = S1}) allows the user to pass the $S$-function $S_1$ to the {\tt Ode1a} command.

\bigskip


%%%%%%%
\noindent
\subsection{Command: {\tt Hfunction}}
\label{hfuncom}

\noindent {\it Feature:} This command tries to find the solution of the associated 1ODE.
\bigskip

\noindent
{\it Calling sequence, parameters and extra parameters:} Same as above.

\bigskip


%%%%%%%
\noindent
\subsection{Command: {\tt PDEsol}}
\label{pdeascom}

\noindent {\it Feature:} This command constructs and tries to solve the 1PDE that relates the $H$-function to the first integral $I$.
\bigskip

\noindent
{\it Calling sequence, parameters and extra parameters:} Same as above.

\bigskip


%%%%%%%
\noindent
\subsection{Command: {\tt Solde}}
\label{invcom}

\noindent {\it Feature:} This command tries to find a Liouvillian general solution for the 1ODE (or, equivalently, first integral of the vector field $\chi$). It returns a $I(x,y,\theta)=c$, where $\theta$ is an elementary function of $(x,y)$.
\bigskip

\noindent
{\it Calling sequence:}

\begin{verbatim}
[> Solde(ode);
\end{verbatim}

\noindent
{\it Parameters:} The same as above.
\medskip

\noindent
{\it Extra parameters:} The same as above and:
\medskip

\hspace
\parindent
{\tt II} - If this parameter is present in the command call, then the command returns a list with the Liouvillian first integral of the vector field $\chi$ and the $z \rightarrow \theta$ transformation: $[I(x,y,z),\,z=\theta(x,y)]$. 
\bigskip

\noindent
{\it Synopsis:}
\smallskip
\smallskip

The command {\tt Solde} is the final part of the process, i.e., it returns (if possible) the general solution of the 1ODE. With the parameter {\tt II} the user can make the command to return the first integral of the vector field $\chi$ (together with $z \rightarrow \theta$ transformation).


\bigskip


\begin{obs}
\label{extrapar}
The extra parameters are not mandatory -- some are even redundant. For example, if we are to provide the $S$-function we do not need to provide the degree of the polynomial candidate $P_c$. If these redundant parameters are supplied simultaneously, the package will take care of doing (hopefully) the best choice.
\end{obs}

\medskip

%%%%%%%%%%%%%%%%
\section{Example of the usage of the package commands}
\label{exusagecom}

In this section we will show the commands of the $LeapsLS$ package simulating its use on a Maple platform to exemplify in a practical way its operation. 

\medskip

\noindent
Consider the 1ODE (example presented in the section \ref{mffi})

\begin{equation}
\label{expackode1}
\frac{dy}{dx}={\frac {{{\rm e}^{x}}{x}^{3}{y}^{2}+{{\rm e}^{x}}{x}^{2}{y}^{2}+2\,{
{\rm e}^{x}}{x}^{2}y+{{\rm e}^{x}}xy+{{\rm e}^{x}}x+{y}^{2}+{{\rm e}^{
x}}}{{x}^{2}{y}^{2}+{{\rm e}^{x}}{x}^{2}+xy+1}}.
\end{equation}
%
After opening a Maple session, we will load the following packages:
%
\begin{verbatim}
[> with(DEtools):  read(`LeapsLS.txt`):
\end{verbatim}
The {\tt :} sign after the command line avoids printing (on-screen) of the result. The {\it DEtools} package loads several commands to handle ODEs. The {\tt read (`LeapsLS.txt`):} command loads our package. We will load 1ODE (\ref{expackode1}) by typing
\begin{verbatim}
[> _1ode := diff(y(x),x) = (exp(x)*x^3*y(x)^2+exp(x)*x^2*y(x)^2
+2*exp(x)*x^2*y(x)+exp(x)*x*y(x)+exp(x)*x+y(x)^2+exp(x))/(x^2*y
(x)^2+exp(x)*x^2+x*y(x)+1):
\end{verbatim}
%
First, we can find the associated vector field:
%
\begin{verbatim}
[> X := Xi(_1ode):
[> fgh := Xi(_1ode,FGH);
\end{verbatim}
%
The command {\tt  X := Xi(\_1ode):} assigns $X$ as the vector field and {\tt fgh := Xi(\_1ode,FGH);} assigns $fgh$ as a list with its coefficients $[f,g,h]$:
%
$$
fgh := [{x}^{2}{y}^{2}+z{x}^{2}+xy+1,z{x}^{3}{y}^{2}+z{x}^{2}{y}^{2}+2\,z{x}^
{2}y+zxy+zx+{y}^{2}+z,z \left( {x}^{2}{y}^{2}+\right.$$
\begin{equation}
\left. z{x}^{2}+xy+1 \right) ] \nonumber
\end{equation}		
%
We can find the $S$-function $S_1$ by typing
\begin{verbatim}
[> S1 := op(Sfunction(_1ode,Deg=5));
\end{verbatim}
\begin{equation}
S1 := -{\frac {{x}^{2}{y}^{2}+{x}^{2}z+xy+1}{x \left( xy+1 \right) ^{2}}}
\end{equation}		
From $S_1$ we can determine the associated 2ODE:
\begin{verbatim}
[> M0N0 := Ode2(_1ode,SF=S1);
\end{verbatim}
%
\begin{equation}
M0N0 := [ \left( zx-y \right)  \left( zx+y \right) ,x \left( xy+1 \right) ^{2}]
\end{equation}		
%
\begin{verbatim}
[> _2odeas := M0N0 [1]/M0N0 [2];
\end{verbatim}
\begin{equation}
\_2odeas := {\frac { \left( zx-y \right)  \left( zx+y \right) }{x \left( xy+1
 \right) ^{2}}}
\end{equation}		
In order to obtain the polynomials $P_1$, $P_2$ and $P_3$, we simply add the parameter {\tt PS} to the entry of the {\tt Ode2} command:
\begin{verbatim}
[> PPP := Ode2(_1ode,SF=S1,PS);
\end{verbatim}
%
\begin{equation}
PPP := [z{x}^{2}{y}^{2}+zxy+{y}^{2}+z,-{x}^{2}{y}^{2}-z{x}^{2}-xy-1,x \left( 
xy+1 \right) ^{2}]
\end{equation}		
%
Using the function $S_1$ we can (by applying the $S$-function method -- see \cite{Noscpc2019}) find the first integral of the vector field $\chi$ (by applying the command {\tt Solde} with the parameter {\tt II}):
\begin{verbatim}
[>  Inv := Solde(_1ode,SF=S1,II);
\end{verbatim}
\begin{equation}
Inv := [{{\rm e}^{ \left( xy+1 \right) ^{-1}}} \left( zx-y \right) ,z={
{\rm e}^{x}}]
\end{equation}
Finally, we can (by doing the substitution $z \rightarrow \theta$ or apllying directly the {\tt Solde} command) obtain the general solution of the 1ODE (\ref{expackode1}):
\begin{verbatim}
[> sol1ode := Solde(_1ode,SF=S1);
\end{verbatim}
%
\begin{equation}
sol1ode := {{\rm e}^{ \left( xy+1 \right) ^{-1}}} \left( {{\rm e}^{x}}x-y
 \right) ={\it \_C1}
\end{equation}
We can now test the solutions found:
\begin{verbatim}
[> DX := Dx(_1ode):
[> X(Inv[1]);
[> DX(lhs(sol1ode));
\end{verbatim}
$$ 0 $$
$$ 0 $$

\subsection{Special features}
\label{sfia}
\hspace\parindent
In this section, we present some special features of the package we use to solve / study more complicated cases. 
As we can see, in the ten 1ODEs shown in the section \ref{hard1odes} we have that $\theta$ is a function of only one variable ($x$ or $y$). In these cases, as shown in the section \ref{api} (theorem \ref{theoAimp}), we can use an improvement that allows us to calculate only one polynomial instead of two ($FastLs$ algorithm). However, due to some technical details, we decided to leave the first part of the $FastLs$ algorithm in the hands of the package user. In this section we will show you how to use the commands in the package to achieve the solution in more difficult cases.

Consider the 1ODE given by:
\begin{equation}
\label{p1ode}
\frac{dy}{dx}=\frac{\left( 3\,{x}^{5}{y}^{2}-6\,\ln  \left( {\frac {x}{y}} \right) {x}^{4
}y+3\, \left( \ln  \left( {\frac {x}{y}} \right)  \right) ^{2}{x}^{3}-
{x}^{4}y+{x}^{3}+x{y}^{2}-y \right) y}
{{x} \left( {x}^{4}y+{x}^{2}{y}^{3}-2\,\ln  \left( {\frac {x}{y}} \right) x{y}^{2}+ \left( \ln 
 \left( {\frac {x}{y}} \right)  \right) ^{2}y+{x}^{3}-x{y}^{2}-y \right) }.
\end{equation}
After loading the package and the 1ODE (\ref{p1ode}), we can try (first) to find the $S$-function using the {\tt Sfunction} command. However, if we type 
\begin{verbatim}
[> t0 := time(): S1 := op(Sfunction(_1ode,Deg=7)); time()-t0;
\end{verbatim}
%
we get
\begin{equation}
S1 := \frac{1}{y}
\end{equation}
$$11.812$$		
%
after $\approx$ 12 seconds (and spending 215 MB of memory). However, if we type 
%
\begin{verbatim}
[> sol := Solde(_1ode,SF=S1);
\end{verbatim}
%
The output is
\begin{equation}
sol := 
\end{equation}		
%
which would probably leave us intrigued, because as the $S$-function is simple the method should have found the solution. What is actually happening is that the following invariant is found for the vector field associated $\chi$: $I = -\ln(x)+\ln(y)+z$. Now, if we pay attention to the fact that the function $\theta$ present in 1ODE is $\theta = \ln(x/y)$, making the substitution $ z \rightarrow \theta $ we will obtain $I = -\ln(x)+\ln(y)+\ln(x/y)=0$ that represents the identity $\theta = \ln(x/y)$ (and therefore, as it is not an answer to what we seek, the program discards it). We can try using the parameter {\tt DegMNP}  which is, actually, a list with three integers {\tt [degM,degN,degP]} which gives to the {\tt Sfunction} command the degrees of the polynomials that are candidates for $M_0$, $N_0$ and $P_2$. This, however, results in a very large consumption of resources and on a memory `overflow'.

However, we can still succeed using the command {\tt Tr1ode} to transform the 1ODE (\ref{p1ode}) into an 1ODE such that $\theta$ is a function of just one variable. We can do this by typing (for example)
\begin{verbatim}
[> _1odetr := Tr1ode(_1ode,TR=[x=x*y,y=y]);
\end{verbatim}
and obtaining a transformed 1ODE given by:
\begin{equation}
\label{1odetra}
\!\!\frac{dy}{dx}\! =\!{\frac {-y \left( 3\,{x}^{5}{y}^{6}-6\,{x}^{4}\ln\left( x \right) {y
}^{4}-{x}^{4}{y}^{4}+3\,{x}^{3} \ln\left( x \right)^{2}{y}^{2}+{x}^{3}{y}^{2}+x{y}^{2}-1 \right) }{3{x}^{6}\!{y}
^{6}\!-\!6{x}^{5}\!\ln\! \left(\! x\! \right)\!{y}^{4}\!-\!2{x}^{5}\!{y}^{4}\!+\!3{x}^{4}\! \ln\!  \left(\! x \!\right)^{\!2}\!{y}^{2}\!-\!{x}^{3}\!{y}^{4}\!+\!2\!\ln \! \left( \!x\! \right)\! x^2\!{y}^{2}\!+\!2x^2\!{y}^{2}\!-\!x\! \ln\!  \left( \!x \!\right)^{\!2}}}
\end{equation}
So, in the transformed 1ODE, we have that $\theta = \ln(x)$. Now, we can apply again the command {\tt Sfunction} (and, in this case, using the $FastLs$ algorithm)
\begin{verbatim}
[> t0:=time(): s1tr := Sfunction(_1odetr[1],Deg=11); time()-t0;
\end{verbatim}
that returns us
\begin{equation}
s1tr:=\left[{\frac {3{x}^{5}\!{y}^{6}\!-\!6{x}^{4}\!z{y}^{4}\!-\!2{x}^{4}\!{y}^{4}\!+\!3{x}^
{3}\!{y}^{2}\!{z}^{2}\!-\!{x}^{2}\!{y}^{4}\!+\!2zx{y}^{2}\!+\!2x{y}^{2}\!-\!{z}^{2}}{y
 \left( {x}^{3}{y}^{2}-1 \right) }}\right]
\end{equation}
$$2.375$$
i.e., the algorithm finds the $S$-function in less than $3$ seconds (and consuming only 30 MB). At this point, we just have to type
\begin{verbatim}
[> soltr := Solde(_1odetr[1],ST=op(s1tr));
\end{verbatim}
to obtain the general solution of the transformed 1ODE (\ref{1odetra}) in $0.3$ seconds:
\begin{equation}
soltr:={\frac {-x{y}^{2}\!\ln\!  \left( {x}^{3}\!{y}^{2}\!-\!1\! \right)\! -\!x{y}^{2}\!\ln\! 
 \left( y \right)\! +\!\ln\!  \left( x \right) \ln\!  \left( {x}^{3}\!{y}^{2}\!-\!1\!
 \right)\! +\!\ln \! \left( x \right) \ln\!  \left( y \right)\! -\!1}{x{y}^{2}-
\ln  \left( x \right) }}={\it \_C1}
\end{equation}
and
\begin{verbatim}
[> Itr := Solde(_1odetr[1],ST=op(s1tr),II);
\end{verbatim}
to obtain the first integral of the transformed vector field:
\begin{equation}
Itr:=\left[{\frac {x{y}^{2}\!\ln\!  \left( {x}^{3}{y}^{2}\!-\!1 \right)\! +\!x{y}^{2}\!\ln\! 
 \left( y \right) \!-z\ln \! \left( {x}^{3}{y}^{2}\!-\!1 \right)\! -z\ln\! 
 \left( y \right) \!+\!1}{-x{y}^{2}+z}},z=\ln\!  \left( x \right) \right]
\end{equation}
The command {\tt Tr1ode} returns (as second output) the inverse transformation. So, by typing 
\begin{verbatim}
[> itr := _1odetr[2]:
[> Ior := subs(itr,Itr[1]);
\end{verbatim}
we obtain the first integral of the vector field associated with the 1ODE (\ref{p1ode}):
\begin{equation}
Ior:= \frac{xy\ln  \left( {\frac {{x}^{3}}{y}}-1 \right) +xy\ln  \left( y
 \right) -z\ln  \left( {\frac {{x}^{3}}{y}}-1 \right) -z\ln  \left( y
 \right) +1}{z-yx}
\end{equation}
We can check this with
\begin{verbatim}
[> X := Xi(_1ode):
[> X(Ior);
\end{verbatim}
$$0$$
Finally, doing the substitution $z \rightarrow \ln(x/y)$, we find the general solution of the original 1ODE (\ref{p1ode}):
\begin{verbatim}
[> DX := Dx(_1ode):
[> sol := subs(z=ln(x/y),Ior=_C1);
[> DX(Ior);
\end{verbatim}
\begin{equation}
sol:=  \frac{xy\ln  \left( {\frac {{x}^{3}}{y}}-1 \right) +xy\ln  \left( y
 \right) -\ln(\frac{x}{y})\ln  \left( {\frac {{x}^{3}}{y}}-1 \right) -\ln(\frac{x}{y})\ln  \left( y
 \right) +1}{\ln(\frac{x}{y})-yx}=\_C1
\end{equation}
$$0$$


\begin{thebibliography}{25}

\bibitem{Liena}
A. Liénard, 
{ Revue générale de l'électricité}, 
{\bf 23}, (1928), 901- 912, 946-954.

\bibitem{LeSm}
N. Levinson and O. Smith, 
{\it A general equation for relaxation oscillations},
Duke Mathematical Journal, {\bf 9}, (1942), 382-403.

\bibitem{SaGaRa}
S. Saha, G. Gangopadhyay, D.S. Ray, 
{\it Systematic designing of bi-rhythmic and tri-rhythmic models in families of Van der Pol and Rayleigh oscillators},
Communications in Nonlinear Science and Numerical Simulation, {\bf 85}, (2020).
https://doi.org/10.1016/j.cnsns.2020.105234

\bibitem{ChSeLa}
V.K. Chandrasekar, M. Senthilvelan and M. Lakshmanan, {\it On the
complete integrability and linearization of certain second order
nonlinear ordinary differential equations}. Proceedings of the Royal
Society London Series A, {\bf 461}, Number 2060, 2005.

\bibitem{Val}
C. Valls,
{\it Liouvillian integrability of some quadratic Liénard polynomial differential systems},
Rendiconti del Circolo Matematico di Palermo Series 2 (2019) 68:499–519
https://doi.org/10.1007/s12215-018-0374-6

\bibitem{Pol}
B. Van der Pol, {\it Sur les oscillations de relaxation},
Philos. Mag. 2 (1926) 978–992.

\bibitem{Pol2}
B. van der Pol, 
{\it On relaxtion-oscillations}, 
The London, Edinburgh and Dublin Philosophical Magazine and Journal of Science, {\bf 2}, (1927) 978-992.

\bibitem{Pol3}
B. Van der Pol, {\it The non linear theory of electrical oscillations},
Proc. Inst. Radio Eng. 22 (1934) 1051–1086.

\bibitem{KuZa}
N.A. Kudryashov, D.I. Sinelshchikov,
{\it On the criteria for integrability of the Liénard equation}, 
Applied Mathematics Letters, {\bf 57}, (2016) 114–120.
https://doi.org/10.1016/j.aml.2016.01.012

\bibitem{GaBrRo}
M.L. Gandarias, M.S. Bruzón, M. Rosa,
{\it Nonlinear self-adjointness and conservation laws for a generalized Fisher equation}, 
Commun. Nonlinear Sci. Numer. Simulat., {\bf 18}, (2013) 1600–1606.
http://dx.doi.org/10.1016/j.cnsns.2012.11.023

\bibitem{ChrL}
C. Christopher,
{\it An algebraic approach to the classification of centers in polynomial Liénard systems}, 
J. Math. Anal. Appl., {\bf 229}, (1998).
https://doi.org/10.1006/jmaa.1998.6175

\bibitem{LlMeTe}
J. Llibre, A.C. Mereu, M.A. Teixeira,
{\it Limit cycles of the generalized polynomial Liénard differential equations}, 
Math. Proc. Cambridge Philos. Soc., {\bf 148}, (2009) 363–383.

\bibitem{PaBiSeLa}
S.N. Pandey, P.S. Bindu, M. Senthilvelan and M. Lakshmanan,
{\it A group theoretical identification of integrable cases of the Liénard-type equation $\ddot{x} + f(x) dot(x) + g(x)=0$. I. Equations having nonmaximal number of Lie point symmetries}, 
J. Math. Phys., {\bf 50}, (2009). 
https://doi.org/10.1063/1.3187783

\bibitem{LlVaL}
J. Llibre, C. Valls,
{\it On the local analytic integrability at the singular point of a class of Liénard analytic differential systems}, 
Proc. Amer. Math. Soc., {\bf 138}, (2010).
https://www.jstor.org/stable/40590614

\bibitem{ZhWa}
L.-H. Zhang and Y. Wang,
{\it A note on periodic solutions of a forced Liénard-type equation}, 
Anziam J., {\bf 51}, (2010).
doi:10.1017/S1446181110000805

\bibitem{GiLl}
J. Giné and J. Llibre,
{\it Weierstrass integrability in Liénard differential systems}, 
J. Math. Anal. Appl., {\bf 377}, (2011) 362–369.
doi:10.1016/j.jmaa.2010.11.005

\bibitem{Gin}
J. Giné,
{\it Liénard equation and its generalizations}, 
Int. J. Bifurc. Chaos Appl. Sci. Eng. {\bf 27} (6) (2017).

\bibitem{LlZh}
J. Llibre and X. Zhang,
{\it Limit cycles of the classical Liénard differential systems: A survey on the Lins Neto, de Melo and Pugh’s conjecture}, 
Expo. Math., {\bf 35}, (2017) 286–299.
http://dx.doi.org/10.1016/j.exmath.2016.12.001

\bibitem{Dem}
M.V. Demina,
{\it Invariant algebraic curves for Liénard dynamical systems revisited}, 
Applied Mathematics Letters, {\bf 84}, (2018).
https://doi.org/10.1016/j.aml.2018.04.013

\bibitem{CaGu}
J.F. Cariñena and P. Guha,
{\it Nonstandard Hamiltonian structures of the Liénard equation and contact geometry}, 
International Journal of Geometric Methods in Modern Physics, {\bf 16}, (2019).
https://doi.org/10.1142/S0219887819400012

\bibitem{Chr} C. Christopher, 
{\it Liouvillian first integrals of second order polynomial differential equations},
Electron. J. Differential Equations, {\bf 49}, (1999), 7 pp. (electronic).

\bibitem{Nosjpa2002-2}  L.G.S. Duarte, S.E.S.Duarte and L.A.C.P. da Mota, 
{\it Analysing the structure of the integrating factors for first-order ordinary differential equations with Liouvillian functions in the solution}, 
Journal of Physics A: Mathematical and General, {\bf 35} (4), (2002).
DOI:10.1088/0305-4470/35/4/312

\bibitem{Dar}
G. Darboux,
{\it M\'emoire sur les \'equations diff\'erentielles alg\'ebriques du premier ordre et du premier
degr\'e (M\'elanges)},
Bull. Sci. Math. 2\`eme s\'erie 2, 60-96, 2, 123-144, 2, 151-200 (1878).

\bibitem{PrSi}
M. Prelle and M. Singer,
Elementary first integral of differential equations.
{\it Trans. Amer. Math. Soc.}, {\bf  279} 215 (1983).

\bibitem{Sin}
M. Singer,
{\it Liouvillian First Integrals},
Trans. Amer. Math. Soc., {\bf  333} 673-688 (1992).

\bibitem{CaLl}
L. Cairó and J. Llibre,
{\it Darboux Integrability for 3D Lotka-Volterra systems},
J. Phys. A: Math. Gen., {\bf 33} 2395-2406 (2000).

\bibitem{ChLlPaZh} C. Christopher, J. Llibre, C. Pantazi and X. Zhang, 
{\it Darboux integrability and invariant algebraic curves for planar polynomial systems}, 
J. Phys. A, {\bf 35}, (2002) 2457–2476. https://doi.org/10.1088/0305-4470/35/10/310

\bibitem{Nosjpa2002-1} L.G.S. Duarte, S.E.S.Duarte and L.A.C.P. da Mota,
{\it A method to tackle first order ordinary differential equations with
Liouvillian functions in the solution}, in J. Phys. A: Math. Gen.,
{\bf 35} 3899-3910 (2002).

\bibitem{Nosjpa2002-2} L.G.S. Duarte, S.E.S.Duarte and L.A.C.P. da Mota,
{\it Analyzing the Structure of the Integrating Factors for First Order
Ordinary Differential Equations with Liouvillian Functions in the
Solution}, J. Phys. A: Math. Gen., {\bf 35} 1001-1006 (2002).

\bibitem{Noscpc2002}  L.G.S. Duarte, S.E.S.Duarte, L.A.C.P. da Mota and J.F.E.
Skea, {\it Extension of the Prelle-Singer Method and a MAPLE
implementation}, Computer Physics Communications, Holanda, v. 144, n. 1, p. 46-62 (2002).

\bibitem{ChGiGiLl} J. Chavarriga, H. Giacomini, J. Giné and J. Llibre, 
{\it Darboux integrability and the inverse integrating factor}, 
J. Differential Eqs., {\bf 194}, (2003) 116–139. https://doi.org/10.1016/S0022-0396(03)00190-6

\bibitem{Nosjcam2005} J. Avellar, L.G.S. Duarte, S.E.S. Duarte, L.A.C.P. da
Mota, {\it Integrating First-Order Differential Equations with
Liouvillian Solutions via Quadratures: a Semi-Algorithmic Method,}
Journal of Computational and Applied Mathematics {\bf 182}, 327-332, (2005).

\bibitem{Noscpc2007} J. Avellar, L.G.S. Duarte, S.E.S.Duarte and L.A.C.P. da Mota,
{\it Determining Liouvillian first integrals for dynamical systems in the plane},
Computer Physics Communications, {\bf 177}, (2007) 584-596. https://doi.org/10.1016/j.cpc.2007.05.014

\bibitem{ChLlPe} C. Christopher, J. Llibre and J.V. Pereira,
{\it Multiplicity of Invariant Algebraic Curves in Polynomial Vector Fields}, 
Pacific Journal of Mathematics, {\bf 229}, (2007) 63-117. https://doi.org/10.2140/pjm.2007.229.63

\bibitem{ChLlPaWa2} C. Christopher, J. Llibre, C. Pantazi and S. Walcher, 
{\it Inverse problems for invariant algebraic curves: Explicit computations}, 
Proc. Roy. Soc. Edinburgh, {\bf 139}, (2009) 287-302.  https://doi.org/10.1017/S0308210507001175

\bibitem{ChLlPaWa3} C. Christopher, J. Llibre, C. Pantazi and S. Walcher, {\it Darboux integrating factors: inverse problems}, 
J. Differential Equations, {\bf 250}, (2011) 1-25. https://doi.org/10.1016/j.jde.2010.10.013

\bibitem{Zha} X. Zhang, 
{\it Liouvillian integrability of polynomial differential systems}, 
Trans. Amer. Math. Soc., {\bf 368}, (2016) 607-620. https://doi.org/10.1090/S0002-9947-2014-06387-3

\bibitem{FeGi} A. Ferragut and H. Giacomini, 
{\it A New Algorithm for Finding Rational First Integrals of Polynomial Vector Fields}, 
Qual. Theory Dyn. Syst., {\bf 9}, (2010) 89–99. https://doi.org/10.1007/s12346-010-0021-x

\bibitem{BoChClWe}  A. Bostan, G. Chèze, T. Cluzeau and J.-A. Weil, 
{\it Efficient algorithms for computing rational first integrals and Darboux polynomials of planar polynomial vector fields}, 
Mathematics of Computation, {\bf 85}, (2016) 1393-1425. https://doi.org/10.1090/mcom/3007

\bibitem{ChCo}  G. Chèze and T. Combot,
{\it Symbolic Computations of First Integrals for Polynomial Vector Fields}, 
Foundations of Computational Mathematics, (2019). https://doi.org/10.1007/s10208-019-09437-9

\bibitem{Nosjde2021} L.G.S. Duarte and L.A.C.P. da Mota,
{\it An efficient method for computing Liouvillian first integrals of planar polynomial vector fields},
Journal of Differential Equations, {\bf 300}, (2021) 356-385. https://doi.org/10.1016/j.jde.2021.07.045





\bibitem{Noscpc2019}
J. Avellar, M.S. Cardoso, L.G.S. Duarte, L.A.C.P. da Mota
{\it Dealing with Rational Second Order Ordinary Differential Equations where both Darboux and Lie Find It Difficult: The S-function Method},  
Computer Physics Communications, {\bf 234}, (2019) 302-314.

\bibitem{astro}
J.-M. Huré, F. Hersant
{\it A new equation for the mid-plane potential of power law disks}
 Volume 467 / No 3 (June I 2007)  Astronomy \& Astrophysica (A\&A), 467 3 (2007) 907-910



\bibitem{Olv}
P.J. Olver,
{\it Applications of Lie Groups to Differential Equations},
Springer-Verlag, (1986).

\bibitem{Ibr} N.H. Ibragimov,
{\it Elementary Lie Group Analysis and Ordinary Differential Equations},
Wiley: Chichester, (1999).

\bibitem{BlAn}
G.W. Bluman and S.C. Anco,
{\it Symmetries and Integration Methods for Differential Equations},
Applied Mathematical Series {\bf Vol. 154}, Springer-Verlag, New York, (2002).

\bibitem{Sch}
F. Schwarz,
{\it Algorithmic Lie Theory for Solving Ordinary Differential Equations},
Chapman  Hall / CRC -- Taylor and Francis Group, (2008).

\bibitem{Noscpc1997} E.S. Cheb-Terrab, L.G.S. Duarte and L.A.C.P. da Mota,
{\it Computer Algebra Solving of First Order ODEs Using Symmetry Methods}.
Comput.Phys.Commun., {\bf 101}, 254, (1997).

\bibitem{Noscpc1998} E.S. Cheb-Terrab, L.G.S. Duarte and L.A.C.P. da Mota,
{\it Computer Algebra Solving of Second Order ODEs Using Symmetry Methods}.
Comput.Phys.Commun., {\bf 108}, 90, (1998).

\bibitem{AbGu}
B. Abraham-Shrauner and A. Guo,
{\it Hidden Symmetries Associated with the Projective
Group of Nonlinear First-Order Ordinary Differential Equations}.
J. Phys. A: Math.Gen., {\bf 25}, 5597-5608, (1992).

\bibitem{GoLe} K.S. Govinder and P.G.L Leach,
{\it A group theoretic approach to a class of second-order ordinary differential
equations not possesing Lie point symmetries}.
J. Phys. A: Math. Gen., {\bf 30}, 2055-68, (1997).

\bibitem{AdMa} A.A. Adam and F.M. Mahomed,
{\it Non-local symmetries of first-order equations}.
IMA J. Appl. Math., {\bf 60}, 187-98, (1998).

\bibitem{GaBrSe} M.S. Bruzón, M.L. Gandarias and M. Senthilvelan,
{\it Nonlocal symmetries of Riccati and Abel chains and their similarity reductions}.
 Journal of Mathematical Physics, {\bf 53}, 023512 (2012).

\bibitem{MuRo}
C. Muriel and J.L. Romero,
{\it New methods of reduction for ordinary differential equations},
IMA J. Appl. Math., {\bf 66}(2), 111-125, (2001).

\bibitem{MuRo2}
C. Muriel and J.L. Romero,
{\it $\,C^{\,\infty}$-Symmetries and reduction of equations without Lie
point symmetries},
J. Lie Theory, {\bf 13}(1), 167-188, (2003).

\bibitem{PuSa}
E. Pucci and G. Saccomandi,
{\it On the reduction methods for ordinary differential equations}.
J. Phys. A: Math. Gen., {\bf 35}, 6145-6155, (2002).

\bibitem{Nuc}
M.C. Nucci,
{\it Jacobi Last Multiplier and Lie Symmetries:
A Novel Application of an Old Relationship},
Journal of Nonlinear Mathematical Physics, {\bf 12}(2), 284-304, (2005).

\bibitem{Nosamc2007}
J. AvellarL.G.S. Duarte, S.E.S.Duarte and L.A.C.P. da Mota,
{\it A semi-algorithm to find elementary first order invariants
of rational second order ordinary differential equations}, Appl. Math. Comp., {\bf
184} 2-11 (2007).

\bibitem{Nosjpa2001}
L.G.S. Duarte, S.E.S.Duarte, L.A.C.P. da Mota and J.F.E. Skea,
{\it Solving second order ordinary differential equations by
extending the Prelle-Singer method}, J. Phys. A: Math.Gen., {\bf
34} 3015-3024 (2001).

\bibitem{LaRa}
M. Lakshmanan and S. Rajasekar, {\it Nonlinear Dynamics:
Integrability, Chaos and Patterns}. New York: Springer-Verlag (2003).


\bibitem{Nosjmp2009}
L.G.S.Duarte and L.A.C.P.da Mota,
{\it Finding Elementary First Integrals for Rational Second Order Ordinary Differential Equations},
J. Math. Phys., {\bf 50}, (2009).

\bibitem{Nosjpa2010}
L.G.S.Duarte and L.A.C.P.da Mota,
{\it 3D polynomial dynamical systems with elementary first integrals},
J. Phys. A: Math. Theor. {\bf 43}, n.6, (2010). doi:10.1088/1751-8113/43/6/065204

\bibitem{GoPiSe}
P. R. Gordoa, A. Pickering and M. Senthilvelan, 
{\it The Prelle-Singer method and Painleve hierarchies}, 
J. Math. Phys., {\bf 55}, 053510 (2014)

\bibitem{TPCSL}
A. K. Tiwari, S. N. Pandey, V. K. Chandrasekar, M. Senthilvelan and M. Lakshmanan,
{\it The inverse problem of a mixed Liénard-type nonlinear oscillator equation from symmetry perspective},
Acta Mechanica, {\bf 227}, Issue 7, (2016) 2039–2051.

\bibitem{MCSL}
R. Mohanasubha, V. K. Chandrasekar, M. Senthilvelan, M. Lakshmanan
{\it Interplay of symmetries and other integrability quantifiers in finite-dimensional integrable nonlinear dynamical systems}
Proc. R. Soc. A., {\bf 472}: 220150847 (2016). http://doi.org/10.1098/rspa.2015.0847

\bibitem{tito1} C Muriel and J L Romero
First integrals, integrating factors and $\lambda$-symmetries of second-order differential equations
Journal of Physics A: Mathematical and Theoretical, {\bf 42}, Number 36

\bibitem{tito2}
C Muriel and J L Romero
Contribution to the Special Issue “Symmetries of Differential Equations: Frames, Invariants and Applications”
Nonlocal Symmetries, Telescopic Vector Fields and $\lambda$-Symmetries of Ordinary Differential Equations
SIGMA 8 (2012), 106

\bibitem{Dav}
J.H. Davenport, Y. Siret and E. Tournier,
{\it Computer Algebra: Systems and Algorithms for Algebraic Computation}.
Academic Press, Great Britain (1993).

\bibitem{Kam}
E. Kamke, 
{\it Differentialgleichungen: L{\"o}sungsmethoden und L{\"o}sungen}, 
Chelsea Publishing Co, New York (1959).

\bibitem{Ran}
Z. Ran, 
{\it One exactly soluble model in isotropic turbulence}, 
Advances and Applications in Fluid Mechanics, {\bf 5}, (2009), 41-47.




\bibitem{PoMa}
B. van der Pol and J. van der Mark, 
{\it The heart beat considered as a relaxation oscillations and an electrical model of the heart}, 
The London, Edinburgh and Dublin Philosophical Magazine and Journal of Science, {\bf 6}, (1928), 763-775.

\bibitem{Fit}
F. Fitzhugh, 
{\it Impulses and physiological states in theoretical models of nerve membranes},
Biophysics Journal, {\bf 1}, (1961), 445-466.

\bibitem{Str}
S. H. Strogatz, 
{\it Nonlinear Dynamics and Chaos}, 
Addison-Wesley, Reading, Massachussets, (1994).

\end{thebibliography}
\end{document}